\documentclass[oneside,english]{amsart}
\usepackage[T1]{fontenc}
\usepackage[latin9]{inputenc}
\usepackage{color}
\usepackage{amssymb}

\makeatletter
\numberwithin{equation}{section} 
\numberwithin{figure}{section} 
  \theoremstyle{plain}
  \newtheorem*{thm*}{Theorem}
  \theoremstyle{plain}
  \newtheorem{thm}{Theorem}[section]
  \theoremstyle{definition}
  \newtheorem{defn}[thm]{Definition}
  \theoremstyle{plain}
  \newtheorem{lem}[thm]{Lemma}
  \theoremstyle{plain}
  \newtheorem{prop}[thm]{Proposition}
  \theoremstyle{remark}
  \newtheorem{rem}[thm]{Remark}
  \theoremstyle{remark}
  \newtheorem*{acknowledgement*}{Acknowledgement}

\@addtoreset{equation}{section}

\usepackage{babel}
\makeatother

\begin{document}

\title[A geometric setting for systems of ordinary differential
equations] {A geometric setting for systems of ordinary
differential equations}

\author[Bucataru]{Ioan Bucataru}
\address{Ioan Bucataru, Oana Constantinescu; Faculty of Mathematics, Al.I.Cuza
University, B-dul Carol 11, Iasi, 700506, Romania}
\urladdr{http://www.math.uaic.ro/\textasciitilde{}bucataru/}
\urladdr{http://www.math.uaic.ro/\textasciitilde{}oanacon/}

\author[Constantinescu]{Oana Constantinescu}

\author[Dahl]{Matias F. Dahl}
\address{Matias F. Dahl; Institute of Mathematics, Aalto
University, FI-00076 Aalto, Helsinki, Finland}
\urladdr{http://www.math.tkk.fi/\textasciitilde{}fdahl/}

\date{\today}

\begin{abstract}
To a system of second order ordinary differential equations (SODE)
one can assign a canonical nonlinear connection that describes the
geometry of the system. In this work we develop a geometric
setting that allows us to assign a canonical nonlinear connection
also to a system of higher order ordinary differential equations
(HODE). For this nonlinear connection we develop its geometry, and
explicitly compute all curvature components of the corresponding
Jacobi endomorphism. Using these curvature components we derive a
Jacobi equation that describes the behavior of nearby geodesics to
a HODE. We motivate the applicability of this nonlinear connection
using examples from the equivalence problem, the inverse problem
of the calculus of variations, and biharmonicity. For example,
using components of the Jacobi endomorphism we express two
Wuenschmann-type invariants that appear in the study of scalar
third or fourth order ordinary differential equations.
\end{abstract}



\maketitle

\section{Introduction}

It is well known that to a system of second order ordinary
differential equations (SODE) one can associate a nonlinear
connection in a canonical way. This nonlinear connection further
induce various geometric objects useful for the study of the SODE
such as the dynamical covariant derivative, parallel transport,
torsions, curvatures, and the Jacobi equation \cite{bucataru07a,
crampin71, grifone72, grifone00, kosambi35, mirona94}. There are
various attempts to define a similar canonical nonlinear
connection also to systems of higher order ordinary differential
equations (HODE), \cite{bucataru00, bucataru05, cantrijn86,
crampin06, deLeon92, miron94, miron97}. A difficulty with this
task is that to a HODE one can associate a number of nonlinear
connections, and each of these offer different information about
the HODE, \cite{andres91, bucataru07b, bm09, byrnes96, catz74,
mironat96}.

The purpose of this work is to develop a unifying geometric
setting for studying systems of second and higher order ordinary
differential equations. Our approach is focused on two ideas of
Kosambi \cite{kosambi35, kosambi36}. First, in \cite{kosambi35},
Kosambi states that a "system of differential equations can be
dealt with geometrically by means of the tensorial operator of
differentiation, called the bi-derivative". In this work we
propose a global expression, called the dynamical covariant
derivative, for the bi-derivative operator introduced by Kosambi.
The dynamical covariant derivative we propose, is associated to a
system of differential equations, of second or higher order, and a
nonlinear connection that is not fixed yet. We determine the
nonlinear connection by requiring the compatibility of the
dynamical covariant derivative with some geometric structure,
introduced by the equations of variation of the system of
differential equations. This follows the statement of Kosambi,
\cite{kosambi36}, that for a system of HODE, "the connection
cannot be determined without the recourse to the equations of
variation".

In Section \ref{sec:sode} we review the geometrical setting for
systems of second order ordinary differential equations. We
describe the canonical nonlinear connection \cite{crampin71,
grifone72}, the dynamical covariant derivative \cite{bd09,
carinena91}, symmetries, newtonoid vector fields \cite{marmo86,
sarlet90}, the Jacobi endomorphism and the Jacobi equation
\cite{crampin96, grifone00}. The purpose of this section is to
present these objects such that their generalization to systems of
higher order ordinary differential equations will be natural.

In Section \ref{sec:hode} we start with a system of HODE,
identified with a vector field $S$, which is called a semispray,
on the $k$-th order tangent bundle $T^kM$ of the configuration
manifold $M$. In this part of the work we show what are the
advantages of using a specific nonlinear connection if we want to
obtain information about the first order variation and the
symmetries of a system of higher order differential equations.
This connection has been proposed in a coordinate form by Miron
and Atanasiu in \cite{mironat96} and has its origins in the work
of Kosambi \cite{kosambi36}. In Proposition \ref{kcharnabla1} and
Theorem \ref{kcharnabla} we fix this nonlinear connection, among
the other ones that can be associated to a semispray of order $k$,
using the following structures: the dynamical covariant derivative
$\nabla$, the Jacobi endomorphism $\Phi$, and an algebraic
structure of the set of newtonoid vector fields of order $k$. In
Theorem \ref{knablasym} we use the chosen nonlinear connection to
characterize the set of symmetries of a semispray of order $k$
using higher order dynamical covariant derivatives and the Jacobi
endomorphism. We provide explicit formulae for all curvature
components of the Jacobi endomorphism. In the last part of this
section we discuss the first order geodesic variation and obtain
the corresponding Jacobi equations for a semispray of order $k$.

In Section \ref{sec:applic} we discuss the applicability of the
theory developed in this paper in various fields, using examples
from the equivalence problem, the inverse problem of the calculus
of variation, and biharmonicity. In Subsection \ref{subsec:prol}
we consider a second order Lagrangian $L_2$ that represents the
second order prolongation of a Riemannian structure \cite{bm09}.
Its variational semispray is a semispray of order $3$ and the
corresponding Jacobi endomorphism encodes information about the
curvature of the Riemannian structure, and their first and second
order covariant derivatives. The second order Lagrangian $L_2$ has
been used also in \cite{cadeo06} to prove that biharmonic curves
are variational. We further show that some invariants that have
been associated to third or fourth order ordinary differential
equations can be expressed in a geometrical way in terms of the
curvature components of the Jacobi endomorphism and their
dynamical covariant derivatives. In Subsection
\ref{subsec:Wuenschmann} we provide a geometric reformulation for
one of these invariants, called the Wuenschmann invariant. The
Wuenschmann invariant represents the obstruction for the existence
of an associated conformal Lorentzian structure on the
$3$-dimensional solution space of a third order differential
equation, \cite{crampin06, frittelli06}. Also, it has been shown
that the Wuenschmann invariant represents an obstruction for an
equivalence problem of third order ordinary differential
equations, \cite{neut02}. In Subsection \ref{subsec:inv} we
provide a geometric reformulation for an invariant that was
proposed by Fels in \cite{fels96} for characterizing variational
fourth order ordinary differential equations. This invariant has
been also considered by Dridi and Neut in \cite{dridi06} for
studying the equivalence problem of fourth order ordinary
differential equations.

\section{Symmetries and first order variation for systems of second order differential
equations} \label{sec:sode} In this section we reformulate some
results about the geometry of a system of second order ordinary
differential equations such that each of them has a generalization
to system of higher order ordinary differential equations.

\subsection{Systems of SODE and geometric structures on $TM$}

Consider $M$ a real, $n$-dimensional manifold, which is a
topological Hausdorff space with countable base that is locally
homeomorphic to $\mathbb{R}^n$, and with $C^{\infty}$-smooth
transition maps. We assume that all objects are smooth where
defined. For a manifold $M$, its tangent bundle $(TM, \pi, M)$ is
denoted by $TM$, the ring of smooth functions on $M$ is denoted by
$C^{\infty}(M)$, and the set of vector fields on $M$ is denoted by
$\mathfrak{X}(M)$. The canonical submersion of the tangent bundle,
$\pi: TM \to M$, induces a natural foliation on $TM$, whose leaves
are tangent spaces $\pi^{-1}(p)=T_pM$, for $p\in M$. Local
coordinates on the base manifold $M$ are denoted by $(x^i)$, while
induced coordinates on $TM$ will be denoted by $(x^i, y^i)$.
Consider also $TM\setminus\{0\} = \{(x,y)\in TM, y\neq 0\}$, the
tangent bundle with the zero section removed.

Throughout the paper we will consider $I$ a nonempty open interval
in $\mathbb{R}$. For a curve $c:I \to M$, $c(t)=(x^i(t))$, denote
by $c': I \to TM$, $c'(t)=(x^i(t), dx^i/dt)$, its tangent lift. We
say that curve $c$ is a \emph{regular curve} if $c'(t)\in
TM\setminus\{0\}$ for all $t\in I$.

The \emph{tangent structure} (or vertical endomorphism) is the
$(1,1)$-type tensor field on $TM$ defined as
$J={\partial}/{\partial y^i}\otimes dx^i$. The \emph{vertical
distribution} is the $n$-dimensional distribution given by  $V:u
\in TM \mapsto V_uTM= \operatorname{Ker} J_u=\operatorname{Im} J_u
= \operatorname{Ker} D_u\pi \subset T_uTM$, where $D_u\pi$ denotes
the tangent map of $\pi$ at $u\in TM$. This distribution is
integrable being tangent to the natural foliation of the tangent
bundle.

In this section we pay attention to the geometry of a system of
second order ordinary differential equations. Such a system of
SODE on $M$ can be represented using a semispray, which is a
special vector field on $TM\setminus\{0\}$.

A \emph{semispray} on $M$ is a vector field $S \in
\mathfrak{X}(TM\setminus\{0\})$ such that any integral curve of
$S$, $\gamma: I \to TM\setminus\{0\}$, is of the form
$\gamma=(\pi\circ \gamma)'$. The reason that a semispray $S$ is
only defined on $TM\setminus\{0\}$ and not on the whole tangent
bundle $TM$ is to include also semisprays induced by Finsler
metrics \cite{bao00, szilasi03}. In induced coordinates $(x,y)$
for $TM\setminus\{0\}$, a semispray $S$ is given by
\begin{eqnarray} S=y^i
\frac{\partial}{\partial x^i} - 2G^i(x,y)\frac{\partial}{\partial
y^i}, \label{semispray}
\end{eqnarray} for some  functions $G^i$ defined on domains of
$(x,y)$.  For an integral curve $\gamma: I \to TM\setminus\{0\}$
of $S$, we say that curve $c:I \to M$, $c=\pi\circ \gamma$ is a
\emph{geodesic} of $S$. Therefore, a regular curve $c:I \to M$ is
a geodesic of $S$ if and only if $S\circ c'=c''$. Locally, a
regular curve $c:I \to M$, $c(t)=(x^i(t))$, is a geodesic of $S$
if and only if it satisfies the system of second order ordinary
differential equations
\begin{eqnarray}
\frac{d^2x^i}{dt^2}+2G^i\left(x, \frac{dx}{dt}\right)=0.
\label{sode}
\end{eqnarray} Thus a semispray describes
systems of SODE with regular curves on $M$ as solutions.

A geometric structure useful for the geometry of a system of SODE
is that of nonlinear connection. An $n$-dimensional distribution
$H: u \in TM\setminus\{0\} \mapsto H_uTM\subset T_uTM$ that is
supplementary to the vertical distribution is called a
\emph{nonlinear connection}, or \emph{horizontal distribution}.
Therefore, for a nonlinear connection, we have the following
decomposition: $T_uTM=H_u\oplus V_u$, for all $u\in
TM\setminus\{0\}$. We will denote by $h$ and $v$ the horizontal
and vertical projectors that correspond to the above
decomposition, and since $h+v=\operatorname{Id}$, either of these
two projectors characterizes the nonlinear connection $H$. We will
denote by $\mathfrak{X}^h(TM\setminus\{0\}) =
h\left(\mathfrak{X}(TM\setminus\{0\})\right)$ and
$\mathfrak{X}^v(TM\setminus\{0\})
=v\left(\mathfrak{X}(TM\setminus\{0\})\right)$ the
$C^{\infty}(TM\setminus\{0\})$-modules of horizontal and
respectively vertical vector fields.

In induced coordinates $(x^i, y^i)$ for $TM$, the vertical
distribution is spanned by ${\partial}/{\partial y^i}$. Therefore,
the horizontal distribution is locally spanned by vector fields of
the form
$$\frac{\delta}{\delta x^i}=h\left(\frac{\partial}{\partial x^i}\right)
=\frac{\partial}{\partial x^i} - N_i^j(x,y)
\frac{\partial}{\partial y^j},$$ for some functions $N^i_j$
defined on domains of induced coordinates on $TM\setminus\{0\}$.
The corresponding horizontal and vertical projectors can be
written as
$$h=\frac{\delta}{\delta x^i} \otimes dx^i, \quad
v=\frac{\partial}{\partial y^i}\otimes \delta y^i. $$ The
$1$-forms $\delta y^i=dy^i+ N^i_jdx^j$ are annihilators for the
horizontal distribution. A change of induced local coordinates
$(x^i, y^i) \to (\widetilde{x}^i(x), \widetilde{y}^i(x,y))$ on
$TM$ induces the following transformation rule for the basis
$\delta/\delta x^i$, $\partial/\partial y^i$ and its dual basis
$dx^i, \delta y^i$
$$\frac{\delta}{\delta x^j}=\frac{\partial \widetilde{x}^i}{\partial
x^j} \frac{\delta}{\delta \widetilde{x}^i}, \quad
\frac{\partial}{\partial y^j}=\frac{\partial
\widetilde{x}^i}{\partial x^j}\frac{\partial}{\partial
\widetilde{y}^i}, \quad d\widetilde{x}^i=\frac{\partial
\widetilde{x}^i}{\partial x^j}dx^j, \quad
\delta\widetilde{y}^i=\frac{\partial \widetilde{x}^i}{\partial
x^j}\delta y^j. $$ Therefore the components of a tensor field on
$TM$, with respect to these bases, will transform as the
components of a tensor field on the base manifold $M$.

Although a semispray $S$ always induces a canonical nonlinear
connection, we will not yet assume any relation between $S$ and
$(h,v)$ at this point. For a semispray $S$ and an arbitrary
nonlinear connection $(h,v)$, consider the vertically valued
$(1,1)$-type tensor field on $TM\setminus\{0\}$,
\begin{eqnarray*} \Phi=-v\circ \mathcal{L}_Sv = v\circ \mathcal{L}_Sh=v\circ
\mathcal{L}_S \circ h,\end{eqnarray*} which will be called the
\emph{Jacobi endomorphism}. Here, for a $(1,1)$-type tensor field
$A$, we denote
$$\mathcal{L}_SA= \mathcal{L}_S\circ A - A\circ \mathcal{L}_S,$$
the Fr\"olicher-Nijenhuis bracket of $S$ and $A$. Locally, the
Jacobi endomorphism can be expressed as follows
\begin{eqnarray} \Phi= R^i_j(x,y) \frac{\partial}{\partial y^i}
\otimes dx^j, \quad R^i_j = 2\frac{\delta G^i}{\delta x^j} -
S(N^i_j) + N^i_kN^k_j. \label{jacobi_end} \end{eqnarray} The
components $R^i_j$ of the Jacobi endomorphism have been considered
be Kosambi \cite{kosambi35} as the second geometric invariant of
the semispray. From this invariant one can determine two other
geometric invariants, namely the curvature of the nonlinear
connection and one component of the curvature of the Berwald
connection \cite{antonelli03}. In the Riemannian context, the
components $R^i_j$ of the Jacobi endomorphism are related to the
curvature $R^i_{jkl}$, \cite{nicolaecsu07}, of the Levi-Civita
connection as follows $R^i_j(x,y)=R^{i}_{kjl}(x)y^ky^l$. A similar
formula holds true also in the Finslerian context,
\cite[(8.15)]{shen01}, where $R^i_j$ is called the Riemann
curvature. For the Jacobi endomorphism $\Phi$ in formula
\eqref{jacobi_end}, we will refer to its components as to the
\emph{curvature components}.

For a pair $(S, (h,v))$, consider the map $\nabla:
\mathfrak{X}(TM\setminus\{0\}) \to
\mathfrak{X}(TM\setminus\{0\})$, given by
\begin{eqnarray} \nabla=h\circ \mathcal{L}_S\circ h + v\circ
\mathcal{L}_S\circ v   =\mathcal{L}_S  + h\circ \mathcal{L}_S h +
v\circ \mathcal{L}_S v \label{nabla}
\end{eqnarray} that will be called the \emph{dynamical covariant
derivative}. By setting $\nabla f=S(f), \textrm{ for } f\in
C^\infty(TM\setminus\{0\})$, using the Leibniz rule, and the
requirement that $\nabla$ commutes with tensor contraction, we
extend the action of $\nabla$ to arbitrary tensor fields and forms
on $TM\setminus\{0\}$, see \cite[Section 3.2]{bd09}. For example,
if $\omega$ is a $1$-form on $TM\setminus\{0\}$, then its
dynamical covariant derivative is given by
\begin{eqnarray}
\left(\nabla \omega\right)(X)=S(\omega(X))-\omega\left(\nabla
X\right). \label{nox} \end{eqnarray} For a $(1,1)$-type tensor
field $A$ on $TM\setminus\{0\}$, its dynamical covariant
derivative is given by
\begin{eqnarray} \nabla A= \nabla\circ A - A\circ \nabla.
\label{nablaa} \end{eqnarray} From first formulae in \eqref{nabla}
and \eqref{nablaa} it follows that $\nabla h= 0$ and $\nabla v=0$.
Hence $\nabla$ preserves the horizontal and vertical distributions
$H$ and $V$. However, $\nabla$ does not act in a similar way on
these two distributions. This can be seen locally as follows.
Using formulae \eqref{nabla} and \eqref{nox} we have
\begin{eqnarray}  \nabla \frac{\delta}{\delta
x^i} =   N^j_i \frac{\delta}{\delta x^j},  & \nabla dx^i  =
-N^i_j dx^j, \label{nabladelta} \\
\nabla \frac{\partial}{\partial y^i}  =  \left(2\frac{\partial
G^j}{\partial y^i}- N^j_i\right) \frac{\partial}{\partial y^j}, &
\nabla \delta y^i  =  -\left(2\displaystyle\frac{\partial
G^i}{\partial y^j}- N^i_j\right)\delta y^j. \nonumber
\end{eqnarray}
The action of the dynamical covariant derivative $\nabla$ on the
components of a horizontal vector field $X=hX$ is given by
\begin{eqnarray} \nabla X=\nabla\left(X^i\frac{\delta}{\delta x^i}\right) =
\nabla X^i \frac{\delta}{\delta x^i}, \quad \nabla
X^i=S(X^i)+N^i_jX^j. \label{nablaxi}\end{eqnarray} In formula
\eqref{nablaxi} both components $X^i$ and $\nabla X^i$ transform,
under a change of induced coordinates on $TM$, as the components
of a vector field from the base manifold $M$. The first order
differential operator $\nabla X^i$ was introduced by Kosambi
\cite{kosambi35} with the name of bi-derivative. The term
"dynamical covariant derivative" was introduced by Cari\~nena and
Martinez in \cite{carinena91}.

Next lemma gives some compatibility conditions between some of the
geometric structures introduced so far. Let us emphasize that in
the next lemma we do not yet assume any relation between the
semispray $S$ and the nonlinear connection $(h,v)$. However, we
will use this lemma to fix a nonlinear connection in Proposition
\ref{charnabla}.

\begin{lem} \label{lem:hlsj} Consider a semispray $S$, a
nonlinear connection $(h,v)$, and the dynamical covariant
derivative $\nabla$ associated to the pair $(S,(h,v))$. Then,
\begin{eqnarray} h\circ \mathcal{L}_S \circ J=-h,  &
J\circ \mathcal{L}_S \circ v=-v, \label{jlshv} \\
\mathcal{L}_SJ + \operatorname{Id} - 2v = \nabla J,  & \nabla
J=2\left(\displaystyle\frac{\partial G^j }{\partial y^i} -
N^j_i\right)\displaystyle\frac{\partial}{\partial y^j}\otimes
dx^i. \label{nablaj}
\end{eqnarray}
\end{lem}
\begin{proof}
For an arbitrary vector field $X\in
\mathfrak{X}(TM\setminus\{0\})$ we have $X+[S,JX]\in
\mathfrak{X}^v(TM\setminus\{0\})$. Composing from the left with
$h$ it follows that $hX=h[JX,S]$, which shows that first formula
\eqref{jlshv} is true. Similarly, composing by $J$, setting
$JX=vZ$ for some $Z\in \mathfrak{X}(TM\setminus\{0\})$ gives the
second equality in formula \eqref{jlshv}.

Using formula \eqref{nabla} and formulae \eqref{jlshv} we obtain
that \begin{eqnarray*} \nabla \circ J=v \circ \mathcal{L}_S \circ
v \circ J = v \circ \mathcal{L}_S \circ J= (\operatorname{Id} - h)
\circ \mathcal{L}_S \circ J = \mathcal{L}_S \circ J +h, \\
J \circ \nabla =J \circ h \circ \mathcal{L}_S \circ h = J \circ
\mathcal{L}_S \circ h = J \circ \mathcal{L}_S \circ
(\operatorname{Id} - v) = J\circ \mathcal{L}_S  +v.\end{eqnarray*}
Using the above two formulae and formula \eqref{nablaa} we obtain
that $\nabla J = \mathcal{L}_SJ + \operatorname{Id} - 2v$, which
is formula \eqref{nablaj}. Second formula \eqref{nablaj} follows
directly from the local expressions \eqref{nabladelta} of the
dynamical covariant derivative.
\end{proof}
We note that from formula \eqref{nablaj}, the $(1,1)$-type tensor
$\nabla J = \mathcal{L}_SJ + \operatorname{Id} - 2v$ represents
the difference between an arbitrary nonlinear connection $(h,v)$,
which is not fixed yet, and the canonical nonlinear connection
associated to a semispray $S$ that will be fixed in Proposition
\ref{charnabla}.

\subsection{Symmetries for systems of SODE}

It is well known that a semispray induces a canonical nonlinear
connection which in turn, determines the five geometric invariants
of the semispray. This is known as KCC theory, after Kosambi,
Cartan, and Chern \cite{antonelli03}. In this section we show that
this connection can be fixed using symmetries of the semispray.
This approach follows the statement of Kosambi \cite{kosambi36}
"the connection cannot be determined without recourse to the
equations of variation". In the next section we show that this
approach of fixing a connection also generalize to systems of
higher order differential equations.

\begin{defn} \label{dynsym} A vector field $X\in \mathfrak{X}(TM\setminus\{0\})$ is a
\emph{dynamical symmetry} of a semispray $S$ if $[S,X]=0$.
\end{defn}
A direct calculation shows that a vector field on
$TM\setminus\{0\}$, locally expressed as
\begin{eqnarray} X=X^i(x,y)\frac{\partial}{\partial x^i} +
Y^i(x,y)\frac{\partial}{\partial y^i}, \label{XTM} \end{eqnarray}
is a dynamical symmetry if and only if  $Y^i=S(X^i)$ and
\begin{eqnarray} S^2(X^i)+X(2G^i)=0. \label{symm1}
\end{eqnarray}
It follows that, while studying dynamical symmetries for a
semispray $S$, the following set of vector fields on
$TM\setminus\{0\}$ plays an important role:
\begin{eqnarray}
\mathfrak{X}^1_S=\left\{X\in \mathfrak{X}(TM\setminus\{0\}), \
X=X^i\frac{\partial}{\partial x^i} +
S(X^i)\frac{\partial}{\partial y^i}\right\}. \label{xs1}
\end{eqnarray} A vector field $X\in \mathfrak{X}^1_S$ is called
a \emph{newtonoid} \cite{marmo86}. Without local coordinates, the
set of newtonoid vector fields can be expressed as follows
\begin{eqnarray} \mathfrak{X}^1_S=\operatorname{Ker}\left(J\circ
\mathcal{L}_S\right)=\operatorname{Im}
\left(\operatorname{Id}+J\circ \mathcal{L}_S\right). \label{xs2}
\end{eqnarray}

For a vector field $X=X^i{\partial}/{\partial x^i}\in
\mathfrak{X}(M)$ the \emph{complete lift}, $X^{1,1}\in
\mathfrak{X}(TM\setminus\{0\})$, is the vector field defined by
\begin{eqnarray} X^{1,1}=X^i(x) \frac{\partial}{\partial x^i} +
\frac{\partial X^i}{\partial x^j}(x)y^j\frac{\partial}{\partial
y^i}.\label{completev}\end{eqnarray} The \emph{vertical lift} of a
vector field $X\in \mathfrak{X}(M)$ is the vector field
$X^{1,0}\in \mathfrak{X}(TM\setminus\{0\})$ given by
$X^{1,0}=J(X^{1,1})$ \cite{yano73}. We denote by
$\mathfrak{X}^{1,1}(TM\setminus\{0\})$, the set of complete lifts
of vector fields on $M$, and by
$\mathfrak{X}^{1,0}(TM\setminus\{0\})=J\left(\mathfrak{X}^{1,1}(TM\setminus\{0\})\right)$,
the set of vertical lifts of vector fields on $M$.

From expression \eqref{xs1} and formula \eqref{completev} it
follows that complete lifts and newtonoid vector fields are
related by
\begin{eqnarray*}
\mathfrak{X}^{1,1}(TM\setminus\{0\}) \subseteq \bigcap_{S
\textrm{\ semispray}} \mathfrak{X}^1_S,
\end{eqnarray*}
and equality holds if and only if $\dim M\geq 2$. If $\dim M=1$
the above inclusion is strict due to the fact that the fibers of
$TM\setminus\{0\}$ are not connected.

\begin{defn} \label{liesym} A vector field $X\in \mathfrak{X}(M)$ is a
\emph{Lie symmetry} of a semispray $S$ if its complete lift,
$X^{1,1}$, is a dynamical symmetry, which means that
$[S,X^{1,1}]=0$.
\end{defn}
Next, we define a $C^{\infty}$-module structure on the set of
newtonoids. This module structure was introduced in
\cite{sarlet90} for discussing adjoint symmetries for systems of
ordinary differential equations. In proposition \ref{charnabla} we
show that this structure can be used to characterize the canonical
nonlinear connection associated to a semispray.

\begin{rem} \label{rem:xs} For $f\in C^{\infty}(TM\setminus\{0\})$ and $X\in
\mathfrak{X}^1_S$, we define the product
$$f\ast X=\left(\operatorname{Id}+J\circ \mathcal{L}_S \right)(fX) = fX+S(f)JX.$$
\begin{itemize} \item[i)]
With respect to the product $\ast$, the set $\mathfrak{X}^1_S$ has
the structure of a $C^{\infty}(TM\setminus\{0\})$-module, and the
set $\mathfrak{X}^{1,1}(TM\setminus\{0\})$ has the structure of a
$C^{\infty}(M)$-module.
\item[ii)] The maps $J: (\mathfrak{X}^1_S, \ast) \to
(\mathfrak{X}^v(TM\setminus\{0\}), \cdot)$ and $h:
(\mathfrak{X}^1_S, \ast) \to (\mathfrak{X}^h(TM\setminus\{0\}),
\cdot)$ are isomorphisms between
$C^{\infty}(TM\setminus\{0\})$-modules. The map $D\pi:
\left(\mathfrak{X}^{1,1}(TM\setminus\{0\}), \ast\right) \to
(\mathfrak{X}(M), \cdot)$ is an isomorphism between
$C^{\infty}(M)$-modules.
\item[iii)] A vector field $X$ on $TM\setminus\{0\}$ is a newtonoid if and only if it can be
expressed as follows
$$ X=X^i(x,y)\ast \frac{\partial}{\partial x^i}.$$  \end{itemize}
\end{rem}

We have seen that a vector field $X$ on $TM\setminus\{0\}$ is a
dynamical symmetry if and only if it is a newtonoid and satisfies
equation \eqref{symm1}. Therefore, a vector field
$X=X^i(x){\partial}/{\partial x^i}$ on $M$ is a Lie symmetry if
and only if components $X^i$ satisfy the system of equations
\eqref{symm1}. Our aim now is to rewrite equations \eqref{symm1}
such that its terms will have a covariant character. Note that
neither $S^2(X^i)$, nor $X(2G^i)$ in formula \eqref{symm1}, has
such a covariant character. For this we evaluate first the
horizontal and vertical components of a vector field, which behave
as the components of a vector field from the base manifold. This
will allow us to characterize newtonoid vector fields using the
dynamical covariant derivative.

\begin{lem} \label{lem:vnablax} Consider a semispray $S$, a nonlinear
connection $(h,v)$, and the dynamical covariant derivative
$\nabla$ associated to the pair $(S, (h,v))$.  A vector field $X$
on $TM\setminus\{0\}$ is a newtonoid if and only if $v(X)=J(\nabla
X),$ which locally is equivalent to
\begin{eqnarray} X=X^i\frac{\delta}{\delta x^i} + \nabla
X^i\frac{\partial}{\partial y^i}, \label{xnabla}\end{eqnarray} for
some functions $X^i$ defined on the domain of induced coordinates
on $TM\setminus\{0\}$, where $\nabla X^i$ is defined in formula
\eqref{nablaxi}.
\end{lem}
\begin{proof}
Since $J\circ \nabla =  J \circ \mathcal{L}_S +v$ it follows that
$J[S,X]=0$ if and only if $v(X)=J(\nabla X)$. The local formula
\eqref{xnabla} follows by formulae \eqref{nablaxi} and
\eqref{xs1}. \end{proof}

\begin{prop} \label{prop:njn} Consider a semispray $S$, a nonlinear
connection $(h,v)$, and the dynamical covariant derivative
$\nabla$ associated to the pair $(S, (h,v))$. A vector field
$X\in \mathfrak{X}(TM\setminus\{0\})$ is a dynamical symmetry if
and only if $X$ is a newtonoid and satisfies
\begin{eqnarray} \nabla (J\nabla X)+\Phi(X)=0. \label{njn}
\end{eqnarray} \end{prop}
\begin{proof}
A vector field $X\in \mathfrak{X}(TM\setminus\{0\})$ is a
dynamical symmetry if and only if $h[S,X]=0$ and $v[S,X]=0$.
Composing by $J$ shows that the first condition is equivalent to
$J[S,X]=0$, which in turns is equivalent to $X\in
\mathfrak{X}^1_S$. If $X\in \mathfrak{X}(TM\setminus\{0\})$ we
have $$v[S,X]=v[S,vX]+v[S,hX]=\nabla(vX)+\Phi(X).$$ Thus $X$ is a
dynamical symmetry if and only if $X\in \mathfrak{X}^1_S$ and
$\nabla(vX)+\Phi(X)=0$ and the result follows by Lemma
\ref{lem:vnablax}. \end{proof}

We note that in the first term of formula \eqref{njn}, $J\nabla X$
only contain horizontal derivatives of $X$, while the outer
$\nabla$ acts on the vertical vector field $J\nabla X$. According
to formula \eqref{nabladelta}, the two actions might be different.
The next proposition shows that the two actions of the dynamical
covariant derivative $\nabla$ coincide if and only if $\nabla J=0$
is satisfied. Moreover, this condition is only satisfied for one
nonlinear connection. In terms of $S$ this unique nonlinear
connection is explicitly given through its vertical projector
$2v=\operatorname{Id}+\mathcal{L}_SJ$, which is the
\emph{canonical nonlinear connection} in \cite{grifone72}. The
advantage of this method of assigning a nonlinear connection to a
semispray is that the method also generalize to systems of higher
order ODEs. This is the topic of Section \ref{sec:hode}.

\begin{prop} \label{charnabla} Consider a semispray $S$, a
nonlinear connection $(h,v)$, and the dynamical covariant
derivative $\nabla$ associated to the pair $(S,(h,v))$. Then, the
following conditions are equivalent:
\begin{itemize} \item[i)] $\nabla$ restricts to a map
$\nabla: \mathfrak{X}^1_S \to \mathfrak{X}^1_S$ that satisfies the
Leibniz rule with respect to the $\ast$ product;
\item[ii)] $\nabla J=0$;
\item[iii)] $\mathcal{L}_SJ+\operatorname{Id}-2v=0$; \item[iv)] $N^i_j={\partial
G^i}/{\partial y^j}$.
\end{itemize} \end{prop}
\begin{proof}
The formulae in \eqref{nablaj} show that conditions ii), iii) and
iv) are equivalent.

To prove implication $ii) \Rightarrow i)$, let $X\in
\mathfrak{X}^1_S$. According to Lemma \ref{lem:vnablax}, we have
that $vX=J\nabla X$. We apply $\nabla$ to both sides of this
equality, use $\nabla v=0$ and $\nabla J=0$, and obtain
$v\left(\nabla X\right)=\nabla vX= \nabla\left(J\nabla
X\right)=J\nabla\left(\nabla X\right),$ which according to Lemma
\ref{lem:vnablax} shows that $\nabla X\in \mathfrak{X}^1_S$.

Let $f\in C^{\infty}(TM\setminus\{0\})$ and $X\in
\mathfrak{X}^1_S$. Since $f\ast X=fX + S(f)JX$ and $\nabla$
satisfies the Leibniz rule for the $\cdot$ product, it follows
that
\begin{eqnarray}
\label{nastfx} \nabla(f\ast X)=S(f)X + f\nabla X + S^2(f)JX+ S(f)\nabla(JX), \\
\nonumber \nabla(f)\ast X + f\ast \nabla X = S(f)X + f\nabla X +
S^2(f)JX+ S(f)J(\nabla X).
\end{eqnarray}
Using the condition $\nabla J=0$ we obtain $\nabla (f\ast
X)=\nabla(f)\ast X + f\ast \nabla X$ .

To prove implication $i) \Rightarrow ii)$, we will show that the
$(1,1)$-type tensor field $\nabla J$ vanishes on the set
$\mathfrak{X}^1_S\cup \mathfrak{X}^v(TM\setminus\{0\})$, which is
a set of generators for $\mathfrak{X}(TM\setminus\{0\})$. By
formulae  \eqref{nabla} and \eqref{nablaa}, $\nabla J$ vanishes on
$\mathfrak{X}^v(TM\setminus\{0\})$. Since $\nabla$ satisfies the
Leibniz rule with respect to the $\ast$ product, it follows from
the two formulae \eqref{nastfx} that $S(f)(\nabla J)X=0$ for an
arbitrary function $f\in C^{\infty}(TM\setminus\{0\})$ and
arbitrary vector field $X\in \mathfrak{X}^1_S$. Therefore, $\nabla
J=0$ on $\mathfrak{X}^1_S$, and this completes the proof of the
proposition.
\end{proof}

For a semispray $S$, we say that the nonlinear connection
characterized by any of the four conditions of Proposition
\ref{charnabla} is the \emph{canonical nonlinear connection} of
$S$. For the remaining of this section we will consider only this
canonical nonlinear connection induced by a semispray $S$. The
next proposition characterizes the dynamical and Lie symmetries of
a semispray $S$ using the dynamical covariant derivative $\nabla$
associated to $S$ and its canonical nonlinear connection.

\begin{prop} \label{nablasym} Consider a semispray $S$. \begin{itemize} \item[i)] A vector field $X\in
\mathfrak{X}(TM\setminus\{0\})$ is a dynamical symmetry of $S$ if
and only if it is a newtonoid and
\begin{eqnarray}
\nabla^2JX + \Phi(X)=0. \label{syminv} \end{eqnarray} \item[ii)]
Locally, a newtonoid vector field $X\in
\mathfrak{X}(TM\setminus\{0\})$ is a dynamical symmetry of $S$ if
and only if its horizontal components satisfy
\begin{eqnarray}
\nabla^2X^i + R^i_jX^j=0, \label{symcoord} \end{eqnarray} where
$\nabla X^i$ is defined in equation \eqref{nablaxi}. \item[iii)] A
vector field $X\in \mathfrak{X}(M)$ is a Lie symmetry of $S$ if
and only if
\begin{eqnarray}
\nabla^2X^{1,0} + \Phi(X^{1,1})=0. \label{liesyminv}
\end{eqnarray}
\end{itemize}
\end{prop}
\begin{proof}
Since we are using the canonical nonlinear connection, we have
that $\nabla J=0$. The first claim follows using Proposition
\ref{prop:njn}. For the second claim we use local formulae
\eqref{nablaxi} and \eqref{jacobi_end}. Therefore, the local
components of the vertical vector field $\nabla^2JX + \Phi(X)$ are
$\nabla^2X^i + R^i_jX^j$. Hence formulae \eqref{syminv} and
\eqref{symcoord} are equivalent. For the third claim consider a
vector field $X\in \mathfrak{X}(M)$. Its complete lift $X^{1,1}$
is a newtonoid for $S$ and formula \eqref{liesyminv} follows from
formula \eqref{syminv} since $JX^{1,1}=X^{1,0}$.
\end{proof}
Equation \eqref{symcoord} represents the covariant reformulation
of equation \eqref{symm1}, both terms $\nabla^2X^i$ and $R^i_jX^j$
behave as components of a vector field from the base manifold.

\subsection{First order variation and Jacobi fields for systems of SODE}

Consider $c: I \to M$ a geodesic of a semispray $S$. A
\emph{geodesic variation} of $c$ is a smooth map $V: I\times
(-\varepsilon, \varepsilon) \to M$, $V=V(t,s)$ such that
\begin{itemize}
\item[i)] $V(t,0)=c(t)$, for all $t$ in $I$,
\item[ii)] $V(t,s)$ is a geodesic of $S$ for all $s$ in $(-\varepsilon,
\varepsilon)$. \end{itemize}
\begin{defn} \label{def:jacobi}
A vector field $\xi:I \to TM$, along a geodesic $c=\pi \circ \xi$
of a semispray $S$, is called a \emph{Jacobi field} of $S$ if
$\xi$ is the variation vector field, $\xi(t)=\left. \partial_s
V(t,s)\right|_{s=0}$, of some geodesic variation $V$.
\end{defn}
Next we characterize the Jacobi fields of a semispray, and hence
the geodesic variation, in terms of a system of differential
equations, using the dynamical covariant derivative and the
curvature components of the Jacobi endomorphism. We define first
the dynamical covariant derivative of a vector field along a
curve.

Consider $\xi:I \to TM$, $\xi(t)=\left(x^i(t), \xi^i(t)\right)$, a
vector field along a regular curve $c=\pi \circ \xi$. Then the
\emph{complete lift} of $\xi$ is the vector field $\xi^{1,1}:I \to
T(TM\setminus\{0\})$ along $c'$, defined as
$\xi^{1,1}(t)=\left(x^i(t), dx^i/dt, \xi^i(t), d\xi^i/dt \right)$.
The \emph{vertical lift} of $\xi$ is the vector field $\xi^{1,0}:I
\to T(TM\setminus\{0\})$ along $c'$, defined as
$\xi^{1,0}(t)=\left(x^i(t), dx^i/dt, 0, \xi^i(t) \right)$. Since
$\xi$ is a vector field along $c$, for each $t\in I$, we can find
a vector field $X\in \mathfrak{X}(M)$ such that $X\circ c=\xi$
near $t$. Then $\xi^{1,1} = X^{1,1}\circ c'$ and $\xi^{1,0} =
X^{1,0}\circ c'$. Assume now that $c$ is a geodesic, which means
that $c'$ is a geodesic of a semispray $S$. Then, we define
$$ \nabla \xi^{1,0}=\left(\nabla X^{1,0}\right)\circ c'.$$ If locally $\xi=(x^i, \xi^i)$, then $\nabla \xi^{1,0}=\nabla \xi^i \partial/\partial
y^i$, where
\begin{eqnarray}\nabla
\xi^i=\frac{d\xi^i}{dt}+ N^i_j\left(x, \frac{dx}{dt}\right)\xi^j.
\label{nablaxi10}\end{eqnarray} Thus $\nabla \xi^{1,0}$ does not
depend on the choice of $X$.
\begin{prop} \label{prop:invjacobi}
Consider a bounded interval $I$. A vector field $\xi:I \to TM$,
along a geodesic $c=\pi\circ \xi$, is a \emph{Jacobi field} if and
only if
\begin{eqnarray}
\nabla^2\xi^{1,0} + \Phi(\xi^{1,1})=0. \label{jacobixi10}
\end{eqnarray} If locally $\xi(t)=\left(x^i(t), \xi^i(t)\right)$, equation \eqref{jacobixi10} is equivalent to
\begin{eqnarray}
 \nabla^2\xi^i +
R^i_j\left(x(t),\frac{dx}{dt}\right)\xi^j=0. \label{eq:invjacobi}
\end{eqnarray}
\end{prop}
\begin{proof}
Let $\xi:I \to TM$, $\xi(t)=\left(x^i(t), \xi^i(t)\right)$ be a
vector field along the geodesic $c=\pi\circ \xi$. Using formula
\eqref{nablaxi10}, we have
\begin{eqnarray}
\label{jacobi11} & & \nabla^2\xi^{1,0} + \Phi(\xi^{1,1}) =
\\ & &
 \left( \frac{d^2\xi^i}{dt^2} + 2 \frac{\partial G^i}{\partial
y^j}\left(x, \frac{dx}{dt}\right)
  \frac{d \xi^j}{dt} + 2 \frac{\partial G^i}{\partial x^j}\left(x,
  \frac{dx}{dt}\right)\xi^j\right) \left. \frac{\partial}{\partial
y^i}\right|_{c'(t)}.  \nonumber \end{eqnarray} If $\xi$ is a
Jacobi field then there is a geodesic variation $V(t,s)$ such that
$V(t,0)=c(t)$,
$$\frac{\partial V^i}{\partial s}(t,0)=\xi^i(t), \textrm{ and }
\frac{\partial^2 V^i}{\partial t^2}+ 2G^i\left(V, \frac{\partial
V}{\partial t}\right)=0, \forall (t,s)\in I \times (-\varepsilon,
\varepsilon).$$  Differentiating the last equation with respect to
$s$ and setting $s=0$ we obtain that the right hand side in
formula \eqref{jacobi11} vanish and equations \eqref{jacobixi10}
and \eqref{eq:invjacobi} follow.

Conversely, if  $\xi:I \to TM$, $\xi(t)=\left(x^i(t),
\xi^i(t)\right)$ is a solution of the Jacobi equations
\eqref{eq:invjacobi} then right hand side in formula
\eqref{jacobi11} vanish. Using the assumption that the interval
$I$ is bounded and \cite[Theorem 4.4]{bd10} it follows that $\xi$
is a variation vector field for some geodesic variation. Therefore
$\xi$ is a Jacobi field.
\end{proof}
Consider $c(t)=(x^i(t))$ a geodesic of a semispray $S$. Its
tangent vector field $c'(t)=(x^i(t), dx^i/dt)$ is a Jacobi field
of $S$. This can be seen either from the fact that the right hand
side in formula \eqref{jacobi11} vanish for $\xi^i=dx^i/dt$ or
from the fact that $c'$ is the variation vector field of the
geodesic variation $V(t,s)=c(t+s)$. In the Riemannian or the
Finslerian context \cite{bao00, docarmo92} we can also see that
the tangent vector field $c'$ of a geodesic $c$ satisfies the
Jacobi equations \eqref{eq:invjacobi} and hence it is a Jacobi
field since the geodesic equations \eqref{sode} can be written as
$\nabla\left({dx^i}/{dt}\right)=0$ and the curvature components
satisfy $R^i_j\left(x, {dx}/{dt}\right){dx^j}/{dt}=0$,
\cite{shen01, szilasi03}.

For a geodesic $c(t)=(x^i(t))$ of a semispray $S$, the vector
field $(x^i(t), tdx^i/dt)$ is a geodesic of $S$ if and only if
$y^j{\partial G^i}/{\partial y^j}=2G^i$. This condition means that
functions $G^i$ are $2$-homogeneous in the fibre coordinates $y$.
In this case $S$ is called a \emph{spray} and this includes
Riemannian and Finslerian cases.

The second implication of Proposition \ref{prop:invjacobi} is
proved in the Riemannian \cite{docarmo92} or Finslerian context
\cite{bao00} using the exponential map. In the general context of
a semispray, due to the lack of homogeneity, this exponential map
cannot be defined and the proof in \cite[Theorem 4.4]{bd10} uses
the flow of the complete lift of a semispray.

In the covariant form of the Jacobi equations
\eqref{eq:invjacobi}, the components $R^i_j$ of the Jacobi
endomorphism contain information about the geodesic behavior and
the stability of the geodesics of a semispray $S$. If the real
part of the eigenvalues of the curvature components $R^i_j$ are
negative then the geodesics will spread out and will be unstable,
while if they are positive the geodesics will be stable
\cite{punzi09}.

\section{Symmetries and first order variation for systems of higher
differential equations} \label{sec:hode}

In this section we will extend the geometric theory for systems of
SODE, developed in the previous sections, to systems of HODE. A
system of $(k+1)$ order ordinary differential equations can be
associated with a semispray $S$ of order $k$, which is a special
vector field on the $k$-th order tangent bundle. For a semispray
of order $k$ there are various nonlinear connections one can
associate to it \cite{andres91, bm09, byrnes96, catz74,
cantrijn86, crampin06}. In the first part of the section we
consider a semispray of order $k$ and an arbitrary nonlinear
connection. For this pair we define the dynamical covariant
derivatives of first and higher order, and the Jacobi
endomorphism. The main contribution is to describe a particular
nonlinear connection that yields covariant expressions for
symmetries and first order variation of the system of HODE. For
this nonlinear connection we explicitly compute all components of
the Jacobi endomorphism and use these components to provide a
geometric expression for the Wuenschmann invariant
\cite{crampin06, dridi06, neut02}.

\subsection{Systems of HODE and geometric structures on $T^kM$}

The framework for studying systems of $(k+1)$ order ordinary
differential equations on a manifold $M$ is the tangent bundle
$T^{k}M=J_{0}^{k}M$ of order $k\geq 1$ \cite{cantrijn86, deleon85,
miron97, tulczyjew76}. This is the $k$-th order jet bundle of
curves $c$ from a neighborhood of $0\in\mathbb{R}$ to $M$. For a
curve $c:I \to M$, $c(t)=(x^i(t))$, denote by $j^kc: I \to T^kM$,
its $k$-th jet lift, defined as
$$j^kc(t)=\left(x^i(t), \frac{1}{1!}\frac{dx^i}{dt}(t), \dots,
\frac{1}{k!}\frac{d^kx^i}{dt^k}(t)\right).$$ If $c$ is a regular
curve, then $j^kc(t) \in T^kM\setminus\{0\}$ for all $t\in I$.
Local coordinates on $M$ are denoted by $(x^{i})$ and the induced
local coordinates on $T^kM$ are denoted by
$(x^{i},y^{(1)i},\cdots,y^{(k)i})$, where
\[
y^{(\alpha)i}(j_{0}^{k}c)=\left.\frac{1}{\alpha!}\frac{d^{\alpha}(x^{i}(c(t))}{dt^{\alpha}}\right\vert
_{t=0},\quad \alpha\in\{1,..,k\} .\] Let also $y^{(0)i}=x^i$.
Denote $T^0M=M$. For each $\alpha \in \{0,1,...,k-1\}$, the
canonical submersion $\pi^k_{\alpha}: T^kM \to T^{\alpha}M$
induces a natural foliation of $T^kM$. Consider also the subbundle
$T^kM\setminus\{0\}=\{(x,y^{(1)},\cdots,y^{(k)})\in T^kM,
y^{(1)}\neq 0\}$. Note that $T^kM\setminus\{0\} =
\left(\pi^k_1\right)^{-1}\left(TM\setminus\{0\}\right)$.

The \emph{tangent structure} (or vertical endomorphism) of order
$k$ is the $(1,1)$-type tensor field on $T^{k}M$ defined as
\begin{equation} J=\frac{\partial}{\partial y^{(1)i}}\otimes
dx^{i}+\frac{\partial}{\partial y^{(2)i}}\otimes
dy^{(1)i}+\cdots+\frac{\partial}{\partial y^{(k)i}}\otimes
dy^{(k-1)i}.\label{eq:J}\end{equation} It has the following
properties: $J^{k+1}=0,\,\,
\operatorname{Im}J^{\alpha}=\operatorname{Ker}J^{k-\alpha+1},
\alpha\in\{1,...,k\}$.

The foliated structure of $T^{k}M$ determines $k$ regular vertical
distributions
\begin{eqnarray*}
V_{\alpha}(u)  =\operatorname{Ker} D_{u}\pi_{\alpha-1}^{k}
=\operatorname{Im} J_u^{\alpha} =\operatorname{Ker}
J_u^{k-\alpha+1}, \textrm{ for } u\in T^kM,  \alpha
\in\{1,...,k\}.\end{eqnarray*} For each $\alpha\in\{1,...,k\},$
the distribution $V_{\alpha}$ is tangent to the fibers of
$\pi_{\alpha-1}^{k}: (x^{i},y^{(1)i},\cdots,y^{(k)i})
\rightarrow(x^{i},y^{(1)i},\cdots,y^{(\alpha-1)i})$, and hence it
is integrable. Note that $\dim V_{\alpha}=(k-\alpha+1)n$,
$\alpha\in\{1,...,k\}$ and $V_{k}(u)\subset
V_{k-1}(u)\subset\cdots\subset V_{1}(u),$ for each $u\in T^kM$.

In this section we pay attention to the geometry of a system of
$(k+1)$ order ordinary differential equations. Such a system of
HODE on $M$ can be represented using a semispray of order $k$,
which is a special vector field on $T^kM\setminus\{0\}$.

A \emph{semispray of order $k$} on $M$ is a vector field $S\in
\mathfrak{X}(T^kM\setminus\{0\})$ such that any integral curve of
$S$, $\gamma:I \to T^kM\setminus\{0\}$ is of the form
$\gamma=j^k(\pi_0^k\circ \gamma)$. In induced coordinates for
$T^kM\setminus\{0\}$, a semispray of order $k$ is given by
\begin{eqnarray}
S=y^{(1)i}\frac{\partial}{\partial
x^{i}}+2y^{(2)i}\frac{\partial}{\partial
y^{(1)i}}+\cdots+ky^{(k)i}\frac{\partial}{\partial
y^{(k-1)i}}-(k+1)G^{i}\frac{\partial}{\partial
y^{(k)i}},\label{ksemispray}\end{eqnarray} for some functions
$G^i$ defined on domains of induced local charts. For an integral
curve $\gamma: I \to T^kM\setminus\{0\}$ of $S$, we say that curve
$c=\pi_0^k\circ \gamma$ is a \emph{geodesic} of $S$. Therefore, a
regular curve $c:I \to M$ is a geodesic of $S$ if and only if
$S\circ j^kc=(j^kc)'$. Locally, a regular curve $c:I \to M$,
$c(t)=(x^i(t))$, is a geodesic of $S$ if and only if it satisfies
the system of $(k+1)$ order ordinary differential equations
\begin{eqnarray}
\frac{1}{(k+1)!}\frac{d^{k+1}x^i}{dt^{k+1}} +
G^{i}\left(x,\frac{dx}{dt},\cdots,\frac{1}{k!}\frac{d^{k}x}{dt^{k}}\right)=0.\label{kode}
\end{eqnarray}
Thus a semispray of order $k$ describes systems of HODE with
regular curves on $M$ as solutions.

Following the line developed in Section \ref{sec:sode} we will
show that the geometry of a semispray of order $k$ can be
determined from a canonical nonlinear connection, associated to
the semispray. An $n$-dimensional distribution $H_0: u \in
T^kM\setminus\{0\} \mapsto H_0(u)\subset T_uT^kM$, that is
supplementary to the vertical distribution $V_1(u)$, is called a
\emph{nonlinear connection}, or \emph{horizontal distribution}.
Therefore, for a nonlinear connection, we have the following
decomposition $T_uT^kM=H_0(u)\oplus V_1(u)$, for all $u\in
T^kM\setminus\{0\}$. For a nonlinear connection $H_0$, consider
also the $n$-dimensional distributions
$H_{\alpha}=J^{\alpha}(H_0)$, for $\alpha \in \{1,...,k-1\}$.
Therefore, for a nonlinear connection, we have the following
decomposition: \begin{eqnarray}
T_uT^kM=\oplus_{\alpha=0}^{k-1}H_{\alpha}(u)\oplus V_k(u), \forall
u\in T^kM\setminus\{0\}.\label{decomptkm}\end{eqnarray} The set of
$k$ distributions  $H_0, H_1,..., H_{k-1}$ is called a
multi-connection in \cite{saunders02}. We will denote by $h_{0},
h_1,...,h_{k-1}, v_k$ the projectors that correspond to
decomposition \eqref{decomptkm}, which characterizes the nonlinear
connection. Sometimes, for uniform notations, we also denote
$h_k=v_k$. Decomposition \eqref{decomptkm} induces also a
decomposition of the vertical subspaces
\begin{eqnarray}
V_{\alpha}(u)=\oplus_{\beta=\alpha}^{k-1}H_{\beta}(u)\oplus
V_k(u), \forall u\in T^kM\setminus\{0\}, \forall \alpha \in
\{1,...,k-1\}. \label{decompva}\end{eqnarray} Decomposition
\eqref{decompva} implies that
$\operatorname{Id}|_{\operatorname{Im} J^{\alpha}}=
\operatorname{Id}|_{\operatorname{Ker} J^{k-\alpha+1}}=
\sum_{\beta=\alpha}^{k} h_{\beta} $, for all $\alpha \in
\{0,...,k\}$. Therefore, for all $\alpha \in \{0,...,k\}$ we have
the following identities
\begin{eqnarray} \label{jalphahv} J^{\alpha}  = \sum_{\beta=\alpha}^{k} h_{\beta}
\circ J^{\alpha}, \quad J^{k-\alpha+1}  =
\sum_{\beta=0}^{\alpha-1} J^{k-\alpha+1}\circ h_{\beta}.
\end{eqnarray}
For $\alpha, \beta \in \{0,...,k\}$, one can reformulate the above
identities as follows
\begin{eqnarray} \label{hbja}
h_{\beta} \circ  J^{\alpha} = \left\{
\begin{array}{cl}
0, & \textrm{ if } \beta\leq \alpha-1, \\
J^{\alpha} \circ h_{\beta-\alpha}, & \textrm{ if } \beta \geq
\alpha,
\end{array} \right.
\end{eqnarray}
\begin{eqnarray} \label{jahb}
J^{\alpha} \circ h_{\beta} = \left\{
\begin{array}{cl}
0, & \textrm{ if } \alpha+\beta \geq k+1, \\
h_{\alpha+\beta}\circ J^{\alpha}, & \textrm{ if } \alpha + \beta
\leq k.
\end{array} \right.
\end{eqnarray}
For $\alpha \in \{0,...,k-1\}$, we denote by
$\mathfrak{X}^{h_{\alpha}}(T^kM\setminus\{0\}) =
h_{\alpha}\left(\mathfrak{X}(T^kM\setminus\{0\})\right)$,
$\mathfrak{X}^{v_{k}}(T^kM\setminus\{0\}) =
v_{k}\left(\mathfrak{X}(T^kM\setminus\{0\})\right)$, the
$C^{\infty}(T^kM\setminus\{0\})$-modules of horizontal and
respectively vertical vector fields.

For a semispray of order $k$, $S$, and a nonlinear connection
$(h_{\alpha}, v_k)$, consider the vertically valued, $(1,1)$-type
tensor field on $T^kM\setminus\{0\}$
\begin{eqnarray} \Phi=-v_k\circ \mathcal{L}_Sv_k=v_k\circ
\sum_{\alpha=0}^{k-1}\mathcal{L}_S h_{\alpha}= v_k\circ
\sum_{\alpha=0}^{k-1}\left(\mathcal{L}_S \circ h_{\alpha}\right)
\label{phik} \end{eqnarray} that we will call the \emph{Jacobi
endomorphism}. Here, for a $(1,1)$-type tensor field $A$ on
$T^kM\setminus \{0\}$, $\mathcal{L}_SA=\mathcal{L}_S\circ A -
A\circ \mathcal{L}_S$ denotes the Fr\"olicher-Nijenhuis bracket of
$S$ and $A$. We will obtain now the local expression of the Jacobi
endomorphism. On the domain of induced local charts on
$T^kM\setminus\{0\}$, consider the vector fields
\begin{eqnarray} \frac{\delta}{\delta
x^i}=h_0\left(\frac{\partial}{\partial x^i}\right)=
\frac{\partial}{\partial x^i}-\sum_{\beta=1}^k N^{j}_{(\beta)i}
\frac{\partial}{\partial y^{(\beta)j}}, \quad \frac{\delta}{\delta
y^{(\alpha)i}}=J^{\alpha} \left(\frac{\delta}{\delta x^i}\right),
\label{abasis}
\end{eqnarray} for $\alpha \in \{1,...,k\}$, which form a basis for $\mathfrak{X}(T^kM\setminus\{0\})$, adapted to the
decomposition \eqref{decomptkm}. Functions $N^{i}_{(\alpha)j}$ are
locally defined on $T^kM\setminus\{0\}$ and are called the
coefficients of the nonlinear connection. The dual basis to the
basis \eqref{abasis} is given by the following locally defined
$1$-forms on $T^kM\setminus\{0\}$:
\begin{eqnarray} dx^i, \ \delta y^{(\alpha)i}=dy^{(\alpha)i}+
\sum_{\beta=1}^{\alpha}M_{(\beta)j}^i dy^{(\alpha-\beta)j}, \
\alpha \in\{1,...,k\}, \label{dbasis} \end{eqnarray} where
\emph{the dual coefficients} $M^{i}_{(\beta)j}$ are given by
\cite{miron97}
\begin{eqnarray}
M^{i}_{(1)j}=N^{i}_{(1)j}, \  M^{i}_{(\alpha)j}=N^{i}_{(\alpha)j}+
\sum_{\beta=1}^{\alpha-1}N^{i}_{(\alpha-\beta)s}M^{s}_{(\beta)j},
\ \alpha \in\{2,...,k\}. \label{nm} \end{eqnarray}  A change of
induced local coordinates $\left(x^i, y^{(\alpha)i}\right) \to
\left(\widetilde{x}^i(x), \widetilde{y}^{(\alpha)i}
\left(x,y^{(\beta)}\right)\right)$ on $T^kM$ induces the following
transformation rule for the basis \eqref{abasis} and its dual
basis \eqref{dbasis}
\begin{eqnarray*}
\frac{\delta}{\delta x^j}=\frac{\partial \widetilde{x}^i}{\partial
x^j} \frac{\delta}{\delta \widetilde{x}^i}, \ \frac{\delta}{\delta
y^{(\alpha)j}}=\frac{\partial \widetilde{x}^i}{\partial x^j}
\frac{\delta}{\delta \widetilde{y}^{(\alpha)i}}, \
\frac{\partial}{\partial y^{(k)j}}=\frac{\partial
\widetilde{x}^i}{\partial x^j} \frac{\partial}{\partial
\widetilde{y}^{(k)i}},
\ \alpha \in \{1,...,k-1\}\\
d\widetilde{x}^i=\frac{\partial \widetilde{x}^i}{\partial
x^j}dx^j, \  \delta\widetilde{y}^{(\alpha)i}=\frac{\partial
\widetilde{x}^i}{\partial x^j}\delta y^{(\alpha)j}, \ \alpha\in
\{1,..,k\}. \end{eqnarray*} Therefore the components of a tensor
field on $T^kM$, with respect to these bases, will transform as
the components of a tensor field on the base manifold $M$.

The projectors $(h_{\alpha}, v_k)$ that correspond to the
decomposition \eqref{decomptkm} can be expressed, with respect to
the bases \eqref{abasis} and \eqref{dbasis} as follows
$$ h_0=\frac{\delta}{\delta x^i} \otimes dx^i,
h_{\alpha}= \frac{\delta}{\delta y^{(\alpha)i}} \otimes \delta
y^{(\alpha)i}, \alpha \in \{1,...,k-1\},
v_k=\frac{\partial}{\partial y^{(k)i}}\otimes \delta y^{(k)i}. $$
With respect to the bases \eqref{abasis} and \eqref{dbasis} the
tangent structure of order $k$ has the following expression:
\begin{eqnarray} J=\frac{\delta}{\delta y^{(1)i}}\otimes
dx^{i}+\frac{\delta}{\delta y^{(2)i}}\otimes \delta
y^{(1)i}+\cdots+\frac{\partial}{\partial y^{(k)i}}\otimes \delta
y^{(k-1)i}.\label{Jadapted}\end{eqnarray} Locally, the Jacobi
endomorphism can be expressed as follows
\begin{eqnarray}
\Phi=\sum_{\alpha=0}^{k-1}R^i_{(\alpha)j}\frac{\partial}{\partial
y^{(k)i}}\otimes \delta y^{(\alpha)j}. \label{localphi}
\end{eqnarray} Under a change of induced local coordinates on $T^kM$,
the $k$ components $R^i_{(\alpha)j}$ of the Jacobi endomorphism
transform as the components of a $(1,1)$-type tensor field on the
base manifold. We will refer to $R^i_{(\alpha)j}$ as to the
\emph{curvature components}. In Subsection \ref{subsec:prol} we
will show that in the Riemannian context, these components are
functions of the Riemannian curvature and its dynamical covariant
derivatives.

We will provide now explicit formulae for the curvature components
$R^i_{(\alpha)j}$ of the Jacobi endomorphism. For a semispray of
order $k$, $S$, and a nonlinear connection with coefficients
$N^i_{(\alpha)j}$ we introduce the following notations:
\begin{eqnarray}
\nonumber I^i_{(2)j} & =& 2 N^i_{(2)j} - S(N^i_{(1)j})
+ N^s_{(1)j}N^i_{(1)s}, \\
I^i_{(\alpha)j} & =& \alpha N^i_{(\alpha)j} - S(N^i_{(\alpha-1)j})
+ N^s_{(1)j}N^i_{(\alpha-1)s} + \sum_{\beta=2}^{\alpha-1}
I^s_{(\beta)j}N^i_{(\alpha-\beta)s}, \label{ialpha}
\end{eqnarray}
for $\alpha \in \{3,...,k\}$. With respect to these notations we
have the following formulae for the Lie bracket of the semispray
$S$ and the vectors of the basis \eqref{abasis}
\begin{eqnarray}
\nonumber \left[S, \frac{\delta}{\delta x^j}\right] & =&
N^i_{(1)j}\frac{\delta}{\delta x^i} +  \sum_{\beta=2}^k
I^i_{(\beta)j}\frac{\delta}{\delta
y^{(\beta-1)i}}  \\
\label{sdelta} & + & \left\{ (k+1)\frac{\delta G^i}{\delta x^j} -
S(N^i_{(k)j}) + N^s_{(1)j}N^i_{(k)s}+
\sum_{\beta=2}^{k}I^s_{(\beta)j}N^i_{(k+1-\beta)s}\right\}
\frac{\partial}{\partial y^{(k)i}}, \\
\nonumber \left[S, \frac{\delta}{\delta y^{(\alpha)j}}\right] & =&
- \alpha \frac{\delta}{\delta y^{(\alpha-1)j}} +
N^i_{(1)j}\frac{\delta}{\delta y^{(\alpha)i}} +
\sum_{\beta=2}^{k-\alpha} I^i_{(\beta)j}\frac{\delta}{\delta
y^{(\alpha+\beta-1)i}}  \\ \nonumber &  + & \left\{
(k+1)\frac{\delta G^i}{\delta y^{(\alpha)j}} - \alpha
N^i_{(k+1-\alpha)j} - S(N^i_{(k-\alpha)j}) +
N^s_{(1)j}N^i_{(k-\alpha)s} \right. \\ & &  + \left.
\sum_{\beta=2}^{k-\alpha}I^s_{(\beta)j}N^i_{(k+1-\alpha-\beta)s}\right\}
\frac{\partial}{\partial y^{(k)i}}, \nonumber
\end{eqnarray} for each $\alpha \in \{1,..., k-1\}$.
From formulae \eqref{sdelta} it follows that $\delta y^{(\alpha)i}
[S, \delta/\delta y^{(\beta)j}]$ behave as $(1,1)$-type tensors on
$M$, for all $\alpha\neq \beta \in \{0,1,...,k\}$. Therefore,
$I^i_{(\alpha)j}$, for $\alpha \in \{2,...,k\}$, behave as
$(1,1)$-type tensors on $M$.

From formulae \eqref{phik} and \eqref{localphi} it follows that
the curvature components of the Jacobi endomorphism can be
computed as follows
\begin{eqnarray*}
R^i_{(\alpha)j}\frac{\partial}{\partial y^{(k)i}}=v_k\left[S,
\frac{\delta}{\delta y^{(\alpha)j}}\right], \quad \alpha\in
\{0,1,...,k-1\}.
\end{eqnarray*} Using the expressions \eqref{sdelta} we
obtain the following formulae for $R^i_{(\alpha)j}$, for $\alpha
\in \{0,1,..., k-1\}$
\begin{eqnarray}
\nonumber R^i_{(\alpha)j} & =&  (k+1)\frac{\delta G^i}{\delta
y^{(\alpha)j}} - \alpha N^i_{(k+1-\alpha)j} - S(N^i_{(k-\alpha)j})
+ N^s_{(1)j}N^i_{(k-\alpha)s} \\ &+&
\sum_{\beta=2}^{k-\alpha}I^s_{(\beta)j}N^i_{(k+1-\alpha-\beta)s}.
\label{localr} \end{eqnarray}

For a $k$-semispray $S$ and a nonlinear connection $(h_{\alpha},
v_k)$, consider the map $\nabla: \mathfrak{X}(T^kM\setminus\{0\})
\to \mathfrak{X}(T^kM\setminus\{0\})$, given by
\begin{eqnarray} \nabla = \sum_{\alpha=0}^{k-1} h_{\alpha}\circ \mathcal{L}_S\circ h_{\alpha} + v_k\circ
\mathcal{L}_S\circ v_k   = \mathcal{L}_S + \sum_{\alpha=0}^{k-1}
h_{\alpha} \circ \mathcal{L}_S h_{\alpha} + v_k\circ
\mathcal{L}_Sv_k. \label{knabla}
\end{eqnarray} We call $\nabla$ the \emph{dynamical covariant
derivative} associated to the pair $(S,(h_{\alpha}, v_k))$. By
setting $\nabla f=S(f)$, for $f\in
C^{\infty}(T^kM\setminus\{0\})$, using the Leibniz rule, and the
requirement that $\nabla$ commutes with tensor contraction,  we
extend the action of $\nabla$ to arbitrary tensor fields and forms
on $T^kM\setminus\{0\}$. Formula \eqref{knabla} implies that
$\nabla h_{\alpha}= 0$ for every $\alpha \in \{0,...,k-1\}$. Hence
$\nabla$ preserves all distributions $H_{\alpha}$. Similarly, we
obtain that $\nabla v_k=0$, hence it preserves the vertical
distribution $V_k$. From formula \eqref{decompva} it follows that
$\nabla$ preserves also all the vertical distributions
$V_{\alpha}$. Dynamical covariant derivative $\nabla$ acts in a
similar way on the distributions $H_{\alpha}$, but has a different
action on the vertical distribution $V_k$. This can be seen
locally as follows. Using formulae \eqref{sdelta} and
\eqref{knabla} we obtain
\begin{eqnarray} \label{knabladelta}
\nabla \frac{\delta}{\delta y^{(\alpha)i}} & = &
h_{\alpha}\left[S, \frac{\delta}{\delta y^{(\alpha)i}}\right]=
N^j_{(1)i} \frac{\delta}{\delta y^{(\alpha)j}}, \quad \alpha
\in\{0,...,k-1\}
\\ \nonumber
\nabla \frac{\partial}{\partial y^{(k)i}} & = & v_k\left[S,
\frac{\partial}{\partial y^{(k)i}}\right]=
\left((k+1)\frac{\partial G^j}{\partial y^{(k)i}}-k
N^j_{(1)i}\right) \frac{\partial}{\partial y^{(k)j}}.
\end{eqnarray}
The action of the dynamical covariant derivative $\nabla$ on a
horizontal vector field $X=h_{\alpha}X$, for $\alpha \in
\{0,1,...,k-1\}$, is given by
\begin{eqnarray} \nabla X=\nabla \left(X^i\frac{\delta}{\delta
y^{(\alpha)i}}\right)= \nabla X^i \frac{\delta}{\delta
y^{(\alpha)i}}, \ \nabla X^i=S(X^i)+N^i_{(1)j}X^j.
\label{nablaxik1} \end{eqnarray} The first order differential
operator $\nabla X^i$ and its higher order iterations were
considered by Kosambi \cite{kosambi36}. For $k=2$, iterated
actions of $\nabla$ where also  considered in \cite{bm09} in order
to fix a nonlinear connection for a second order Lagrange space.

Next lemma gives some compatibility conditions between the
geometric structures introduced so far. In this lemma we do not
assume any relation between the semispray and the nonlinear
connection. However, we will use this lemma to fix partially the
nonlinear connection later on.
\begin{lem} \label{lem:khlsj} Consider a semispray $S$, a
nonlinear connection $(h_{\alpha},v_k)$, and the dynamical
covariant derivative $\nabla$ associated to the pair
$(S,(h_{\alpha},v_k))$. Consider $\alpha, \beta \in \{0,...,k\}$.
Then, the following formulae are true:
\begin{eqnarray} J \circ \mathcal{L}_SJ^{\alpha}=-\alpha J^{\alpha},
 \label{jbeta} \end{eqnarray}
\begin{eqnarray} \label{halsj}
h_{\beta} \circ \mathcal{L}_S\circ J^{\alpha} = \left\{
\begin{array}{cl}
-(\beta+1)h_{\beta}\circ
J^{\beta}, & \textrm{ if } \beta=\alpha-1 \\
0, & \textrm{ if } \beta<\alpha-1, \end{array} \right.
\end{eqnarray}
\begin{eqnarray} \label{jblsha}
J^{\beta} \circ \mathcal{L}_S\circ h_{\alpha} = \left\{
\begin{array}{cl}
-\alpha J^{k-\alpha} \circ h_{\alpha}, & \textrm{ if } \alpha=k-\beta+1 \\
0, & \textrm{ if } \alpha>k-\beta+1, \end{array} \right.
\end{eqnarray}
\begin{eqnarray} \nabla
J^{\alpha} \circ J^{k-\alpha} = \nabla J^k=\mathcal{L}_SJ^k+
(kh_{k-1}-v_{k})\circ J^{k-1}, \textrm{ if } \alpha \geq 1,
\label{nablajk}
\end{eqnarray}
\begin{eqnarray} \nabla J^{\alpha} = (k+1)\left(\frac{\partial
G^i}{\partial y^{(k)j}} -
N^i_{(1)j}\right)\frac{\partial}{\partial y^{(k)i}}\otimes \delta
y^{(k-\alpha)j}, \textrm{ if } \alpha \geq 1, \label{nablajalpha}
\end{eqnarray}
\begin{eqnarray} \label{lsjk}
 \mathcal{L}_SJ & + & \operatorname{Id} -(k+1)v_k  =  \\
\displaystyle \nabla J & + & \sum_{\alpha=1}^{k-1}
(k+1)\left(\frac{\delta G^j}{\delta y^{(\alpha)i}}-
N^j_{(k+1-\alpha)i}\right)\frac{\partial}{\partial
y^{(k)j}}\otimes \delta y^{(\alpha-1)i}.  \nonumber
\end{eqnarray}
\end{lem}
\begin{proof}
From formula \eqref{eq:J}, a local computation shows that for
every $\alpha \in \{1,...,k\}$ and for every $X \in
\mathfrak{X}(T^kM\setminus\{0\})$ we have
\begin{eqnarray} [S,J^{\alpha}X] - J^{\alpha}[S,X] + \alpha
J^{\alpha-1}X\in \operatorname{Ker}J =
\operatorname{Im}J^k\label{sjbetax1}.
\end{eqnarray} For $\alpha \in \{1,...,k\}$, formula \eqref{jbeta} follows by composing in formula
\eqref{sjbetax1} to the left with $J$. For $\alpha=0$, formula
\eqref{jbeta} follows since $\mathcal{L}_S \operatorname{Id}=0$.

Formulae \eqref{halsj} follow from formula \eqref{sjbetax1} by
composing to the left with $h_{\beta}$ and using the identities
\eqref{hbja}.

Formulae \eqref{jblsha} follow from formula \eqref{sjbetax1} by
composing to the right with $h_{\beta}$ and using the identities
\eqref{jahb}.

From formulae \eqref{knabla}, \eqref{hbja}, \eqref{jahb},
\eqref{halsj}, and \eqref{jblsha} for $\beta=k$, we obtain
\begin{eqnarray*}
\nabla \circ J^k &=& v_{k}\circ \mathcal{L}_S \circ J^k =
\mathcal{L}_S \circ J^k -(h_0+\cdots + h_{k-1}) \circ
\mathcal{L}_S \circ J^k \\ &=& \mathcal{L}_S \circ J^k + k h_{k-1}
\circ J^{k-1}, \\
J^k \circ  \nabla & = & J^{k}\circ \mathcal{L}_S \circ h_0 =
J^{k}\circ \mathcal{L}_S - J^{k}\circ \mathcal{L}_S \circ
(h_1+\cdots + h_{k-1}+v_k) \\ &=& J^{k}\circ \mathcal{L}_S +
J^{k-1} \circ h_1 = J^{k}\circ \mathcal{L}_S + v_k \circ J^{k-1}.
 \end{eqnarray*}
Now the second equality in formula \eqref{nablajk} follows using
the fact that $\nabla J^k= \nabla \circ J^k - J^k \circ \nabla$.

We will prove now the first equality of formula \eqref{nablajk}.
Since $\nabla J^{\alpha}\circ J^{k-\alpha}=\nabla \circ J^k -
J^{\alpha} \circ \nabla \circ J^{k-\alpha}$, it remains to show
that $J^{\alpha} \circ \nabla \circ J^{k-\alpha}=J^k \circ \nabla
$, for $\alpha \in \{1,...,k\}$. Using formula \eqref{knabla} we
have
\begin{eqnarray*} J^{\alpha} \circ \nabla \circ J^{k-\alpha} = J^{\alpha}
\circ h_{k-\alpha} \circ \mathcal{L}_S \circ h_{k-\alpha} \circ
J^{k-\alpha} = v_k \circ J^{\alpha} \circ \mathcal{L}_S \circ
J^{k-\alpha} \circ h_{0}. \label{jnjk}\end{eqnarray*} In formula
\eqref{jbeta}, we replace $\alpha$ by $k-\alpha$, compose to the
left with $J^{\alpha-1}$ and obtain
$$J^{\alpha} \circ \mathcal{L}_S \circ J^{k-\alpha}= J^{k} \circ
\mathcal{L}_S- (k-\alpha)J^{k-1}. $$ Using the above two formulae
and the fact that $v_k\circ J^{k-1}\circ h_0=0$, we obtain
$$ J^{\alpha} \circ \nabla \circ J^{k-\alpha} = J^k \circ \mathcal{L}_S
\circ h_0 = J^k \circ \nabla.$$

Formula \eqref{nablajalpha} follows from expression
\eqref{Jadapted} and the local expression of the dynamical
covariant derivative \eqref{knabladelta}.

By direct calculation it follows that for $\alpha\in\{0,...,k-1\}$
\begin{eqnarray*}
\left(\operatorname{Id} + \mathcal{L}_SJ\right)
\left(\frac{\delta}{\delta y^{(\alpha)i}}\right) &=& (k+1)
\left(\frac{\delta G^j}{\delta y^{(\alpha+1)i}} -
N^j_{(k-\alpha)i}\right)\frac{\partial}{\partial y^{(k)j}},  \\
\left(\mathcal{L}_SJ\right)\left(\frac{\partial}{\partial
y^{(k)i}}\right) & =&  k\frac{\partial}{\partial y^{(k)i}}.
\end{eqnarray*} Using these formulae and \eqref{nablajalpha} for
$\alpha=1$ it follows that formula \eqref{lsjk} is true.
\end{proof}
When $\beta=1$ formulae \eqref{halsj} and \eqref{jblsha} read as
follows
\begin{eqnarray}
 h_{0}\circ \mathcal{L}_S \circ J=-h_{0}, \quad J\circ
\mathcal{L}_S \circ v_k=-kv_k. \label{jlshvk}
\end{eqnarray}
Formula \eqref{lsjk} generalizes formula \eqref{nablaj}. In the
case $k=1$ the tensor $\nabla J$ uniquely determines the tensor $
\mathcal{L}_SJ +\operatorname{Id} -2v$ that fixes the canonical
nonlinear connection. For the case $k>1$, a nonlinear connection
can be fixed using the (1,1)-type tensor $ \mathcal{L}_SJ
+\operatorname{Id} -(k+1)v_k $. Using formula \eqref{lsjk}, in
Theorem \ref{kcharnabla}, we will determine this tensor by fixing
$\nabla J$ and
\begin{eqnarray} \sum_{\alpha=1}^{k-1}(k+1)\left(\frac{\delta G^j}{\delta
y^{(\alpha)i}}- N^j_{(k+1-\alpha)i}\right)\frac{\partial}{\partial
y^{(k)j}}\otimes \delta y^{(\alpha-1)i}.\label{ijphi}
\end{eqnarray} Therefore, in Theorem \ref{kcharnabla}, we will determine
a nonlinear connection for a semispray $S$ by fixing $\nabla J$
and the $(1,1)$-type tensor whose local expression is given in
formula \eqref{ijphi}.

\subsection{Symmetries for systems of HODE} \label{subsec:symhode}

For a semispray of order $k$ there are various nonlinear
connections one can associate to it \cite{andres91, bucataru05,
catz74, cantrijn86, mironat96}. In this section we will show what
are the advantages of using the nonlinear connection proposed by
Miron and Atanasiu \cite{mironat96}. We show that this connection
is very useful for studying the symmetries and the first order
variation of a system of higher order ordinary differential
equations. Moreover, we show that with respect to this nonlinear
connection one can explicitly write down all components of the
Jacobi endomorphism. We use these components to provide a
geometric expression for the Wuenschmann invariant \cite{
crampin06, dridi06, neut02} as well as for an invariant introduced
by Fels \cite{fels96} for the inverse problem of a fourth order
ODE.

\begin{defn} \label{defn:kdynsym} A vector field $X\in \mathfrak{X}(T^kM\setminus\{0\})$ is a
\emph{dynamical symmetry} of a $k$-semispray $S$ if $[S,X]=0$.
\end{defn}

A direct calculation shows that a vector field on
$T^kM\setminus\{0\}$,
\begin{eqnarray} X=X^i\left(x,y^{(\beta)}\right)\frac{\partial}{\partial x^i} +
\sum_{\alpha=1}^k Y^{(\alpha)i}\left(x,y^{(\beta)}\right)
\frac{\partial}{\partial y^{(\alpha)i}}, \label{XTKM}
\end{eqnarray} is a dynamical symmetry if and only if its
components satisfy $\alpha! Y^{(\alpha)i}=S^{\alpha}(X^i)$, for
all $\alpha \in \{1,...,k\}$ and
\begin{eqnarray}
S^{k+1}(X^i)+(k+1)!X(G^i)=0. \label{symmk}
\end{eqnarray}
It follows that, while studying dynamical symmetries for a
$k$-semispray $S$, the following set of vector fields on
$T^kM\setminus\{0\}$ plays an important role:
\begin{eqnarray}
\mathfrak{X}^k_S=\left\{X\in \mathfrak{X}(T^kM\setminus\{0\}),
X=X^i\frac{\partial}{\partial x^i} + \sum_{\alpha=1}^k
\frac{1}{\alpha!}S^{\alpha}(X^i)\frac{\partial}{\partial
y^{(\alpha)i}} \right\}. \label{xsk}
\end{eqnarray} A vector field $X\in \mathfrak{X}^k_S$ is called
a \emph{newtonoid} of order $k$. Next, we will provide global
characterizations, in terms of some differential operators for the
set of newtonoid vector fields $\mathfrak{X}^k_S$ and their images
$J^{\alpha}\left(\mathfrak{X}^k_S\right)$, for $\alpha \in
\{1,...,k\}$.

Consider $\pi^k_S: \mathfrak{X}(T^kM\setminus\{0\}) \to
\mathfrak{X}(T^kM\setminus\{0\})$ the differential operator of
order $k$ defined as
$$\pi^k_S=\operatorname{Id} + \frac{1}{1!}J \circ \mathcal{L}_S +
\cdots + \frac{1}{k!}J^k \circ \mathcal{L}_S^k.
$$
The form of the differential operator $\pi^k_S$ was inspired by
the generalized Cartan operator acting on $1$-forms on $T^kM$,
considered in \cite{cantrijn86, sarlet90}.
\begin{lem} \label{lem33}
For each $\alpha \in \{0,...,k\}$, the set
$J^{\alpha}\left(\mathfrak{X}^k_S\right)$ can be expressed,
without local coordinates, as follows
\begin{eqnarray} J^{\alpha}\left(\mathfrak{X}^k_S\right) =
\operatorname{Ker}\left(J\circ \mathcal{L}_S + \alpha
\operatorname{Id} \right) \cap \operatorname{Im} J^{\alpha}
=\operatorname{Im}(J^{\alpha}\circ \pi^k_S). \label{xsk1}
\end{eqnarray}
\end{lem}
\begin{proof}
First, we prove the two equalities in formula \eqref{xsk1} for
$\alpha=0$. Consider a vector field $X$ on $T^kM\setminus\{0\}$,
given locally by formula \eqref{XTKM}. It follows that $J[S,X]=0$
if and only if its components satisfy $\beta!
Y^{(\beta)i}=S^{\beta}(X^i)$, for all $\beta \in \{1,...,k\}$.
This proves first equality in formula \eqref{xsk1}.

From formula \eqref{jbeta} we have $J\circ \mathcal{L}_S \circ J =
J^2 \circ \mathcal{L}_S -J$. Using this formula recurrently, we
obtain that $(J\circ \mathcal{L}_S) \circ \pi^k_S=0$. Therefore,
$\operatorname{Im} \pi^k_S \subset \operatorname{Ker}\left(J\circ
\mathcal{L}_S\right)$.

For the converse inclusion, consider a vector field $X$ on
$T^kM\setminus\{0\}$ such that $J[S, X]=0$. It follows that
$\mathcal{L}_S^{\beta} X \in \operatorname{Ker} J^{\beta}$, for
all $\beta \in \{1,...,k\}$. Hence, $\left(J^{\beta}\circ
\mathcal{L}_S^{\beta}\right) X=0$ and $X=\pi^k_S(X)$.

Consider now $\alpha \in \{1,...,k\}$. Using the above
considerations it follows that
$J^{\alpha}\left(\mathfrak{X}^k_S\right)=\operatorname{Im}(J^{\alpha}\circ
\pi^k_S)$. We prove now the first equality in formula
\eqref{xsk1}. Consider $Z=J^{\alpha}X$ for $X \in
\mathfrak{X}^k_S$. Therefore $J[S,X]=0$ and using formula
\eqref{jbeta} we have $J[S, J^{\alpha}X] +\alpha J^{\alpha}X=
J^{\alpha+1}[S,X]=0$. It follows that $Z\in
\operatorname{Ker}\left(J\circ \mathcal{L}_S + \alpha
\operatorname{Id} \right) \cap \operatorname{Im} J^{\alpha}$.

For the converse inclusion, consider $Z=J^{\alpha}X$ for $X \in
\mathfrak{X}(T^kM\setminus\{0\})$ such that $J[S,Z]+\alpha Z=0$.
Using again formula \eqref{jbeta} it follows that
$J^{\alpha+1}[S,X]=0$. If $X$ is given locally by formula
\eqref{XTKM}, then its components satisfy $\beta!
Y^{(\beta)i}=S^{\beta}(X^i)$, for all $\beta \in
\{1,...,k-\alpha\}$. It follows that $Z=J^{\alpha}\pi^k_S X \in
J^{\alpha}\left(\mathfrak{X}^k_S\right)$.
\end{proof}

For a vector field $X=X^i\partial/\partial x^i \in
\mathfrak{X}(M)$, the \emph{complete lift}, $X^{k,k}\in
\mathfrak{X}(T^kM\setminus\{0\})$, is the vector field defined as
\begin{eqnarray} X^{k,k}=X^i(x)\frac{\partial}{\partial x^i} +
\sum_{\alpha=1}^k
\frac{1}{\alpha!}d_T^{\alpha}(X^i(x))\frac{\partial}{\partial
y^{(\alpha)i}}, \label{xcomplete}\end{eqnarray} where $d_T$ is the
\emph{Tulczyjev operator} \cite{tulczyjew76}
$$d_T=y^{(1)i}\frac{\partial}{\partial
x^{i}}+2y^{(2)i}\frac{\partial}{\partial
y^{(1)i}}+\cdots+ky^{(k)i}\frac{\partial}{\partial y^{(k-1)i}}.$$
For each $\alpha \in \{0,...,k-1\}$, the \emph{vertical lift} of a
vector field $X\in \mathfrak{X}(M)$ is the vector field
$X^{k,\alpha}=J^{k-\alpha}\left(X^{k,k}\right)\in
\mathfrak{X}(T^kM\setminus\{0\})$, \cite{cantrijn86, deleon95}. We
denote by $\mathfrak{X}^{k,k}(T^kM\setminus\{0\})$, the set of
complete lifts of vector fields on M and by
$\mathfrak{X}^{k,\alpha}(T^kM\setminus\{0\})$, for $\alpha
\in\{0,..., k-1\}$, the set of vertical lifts of vector fields on
M. From expressions \eqref{xsk} and \eqref{xcomplete} it follows
that complete lifts, vertical lifts, and newtonoid vector fields
are related by
\begin{eqnarray} \mathfrak{X}^{k,\alpha}(T^kM\setminus\{0\}) \subseteq  \bigcap_{ S
\textrm{ semispray}} J^{k-\alpha}\left(\mathfrak{X}^k_S\right),
\label{inclusionk}
\end{eqnarray}
for all $\alpha \in \{0,...,k\}$, and the equality holds if and
only if $\dim M\geq 2$ and $\alpha=k$. If $\dim M=1$ the above
inclusion is strict due to the fact that the fibers of
$T^kM\setminus\{0\}$ are not connected.
\begin{defn} \label{defn:kliesym} A vector field $X\in \mathfrak{X}(M)$ is a
\emph{Lie symmetry}, or a point symmetry, of a $k$-semispray $S$
if the complete lift $X^{k,k}$ is a dynamical symmetry, which
means $[S,X^{k,k}]=0$.
\end{defn}
Next, we define a $C^{\infty}$-module structure on the set of
newtonoid vector fields of order $k$. This structure is inspired
by a module structure on a set of $1$-forms associated to adjoint
symmetries of a system of HODE, \cite{sarlet90}. In proposition
\ref{kcharnabla1} we show that this structure can be used to fix a
part of a canonical nonlinear connection one can associate to a
semispray.

\begin{rem} For $f\in C^{\infty}(T^kM\setminus\{0\})$ and $X\in \mathfrak{X}^k_S$, we define
\begin{eqnarray} f\ast X= \pi^k_S(fX) = \sum_{\alpha=0}^k \frac{1}{\alpha!}
S^{\alpha}(f)J^{\alpha}X. \label{kast}
\end{eqnarray}
\begin{itemize} \item[i)] With respect to the product $\ast$, the
set $J^{\alpha}\left(\mathfrak{X}^k_S\right)$ has the structure of
a $C^{\infty}(T^kM\setminus\{0\})$-module and the set
$\mathfrak{X}^{k,\alpha}(T^kM\setminus\{0\})$ has the structure of
a $C^{\infty}(M)$-module, for all $\alpha \in\{0,1,...,k\}$.
\item[ii)] The map $h_{\alpha}
: \left(J^{\alpha}\left(\mathfrak{X}^k_S\right), \ast\right) \to
(\mathfrak{X}^{h_{\alpha}}(T^kM\setminus\{0\}), \cdot)$  is an
isomorphism of $C^{\infty}(T^kM\setminus\{0\})$-modules, for all
$\alpha \in \{0,...,k\}$. The map $D\pi^k_0:
(\mathfrak{X}^{k,k}(T^kM\setminus\{0\}), \ast) \to
(\mathfrak{X}(M), \cdot)$  is an isomorphism of
$C^{\infty}(M)$-modules.
\item[iii)] A vector field $X$ on $T^kM\setminus\{0\}$ belongs
to the set $J^{\alpha}\left(\mathfrak{X}^k_S\right)$ if and only
if it can be expressed as
$$ X=X^i\left(x,y^{(1)},...,y^{(k)}\right)\ast \frac{\partial}{\partial y^{(\alpha)i}},
\quad \forall \alpha \in \{0,1,...,k\}.$$
\end{itemize}
\end{rem}
As we have seen, a vector field $X$ on $T^kM\setminus\{0\}$ is a
dynamical symmetry if and only if it is a newtonoid and satisfies
equations \eqref{symmk}. Our aim now is to rewrite equations
\eqref{symmk} such that its terms have a covariant character. Note
that neither one of the two terms of equation \eqref{symmk},
$S^{k+1}(X^i)$ and $X(G^i)$, has such a covariant character. To
accomplish this goal, we first characterize newtonoid vector
fields in terms of their expressions in the adapted basis
\eqref{abasis}. The components of a newtonoid vector field in the
adapted basis \eqref{abasis} will behave as the components of a
vector field on the base manifold $M$.

\begin{lem} \label{lem:kvnablax} Consider a semispray $S$, a
nonlinear connection $(h_{\alpha},v_k)$, and the dynamical
covariant derivative $\nabla$ associated to the pair
$(S,(h_{\alpha},v_k))$. A vector field $X$ on $T^kM\setminus\{0\}$
is a newtonoid if and only if
\begin{equation} X=X^{i}\frac{\delta}{\delta x^{i}}+ \frac{1}{1!} \nabla^{(1)}
X^{i}\frac{\delta}{\delta y^{(1)i}} +
\frac{1}{2!}\nabla^{(2)}X^{i}\frac{\delta}{\delta y^{(2)i}} +
\cdots + \frac{1}{k!}\nabla^{(k)}X^{i}\frac{\partial}{\partial
y^{(k)i}},\label{eq:symadaptedbasis}\end{equation} where for
$\alpha \in \{1,...,k\}$ we denote
\begin{eqnarray}
\label{eq:nablaiterated} \frac{1}{\alpha !}\nabla^{(\alpha)}X^{i}
& = & \frac{1}{\alpha !} S^{\alpha}(X^{i}) +
\sum_{\beta=1}^{\alpha}\frac{1}{(\alpha-\beta)!}M_{(\beta)j}^{i}S^{\alpha-\beta}(X^{j}).
\end{eqnarray}
For a vector field $X\in J^{k-1}\mathfrak{X}(T^kM\setminus\{0\})$
we have that $X \in J^{k-1}\left(\mathfrak{X}^k_S\right)$ if and
only if
$$ v_k(X)=J(\nabla X).$$ \end{lem}
\begin{proof}
Consider a vector field $X \in \mathfrak{X}(T^kM\setminus\{0\})$.
In view of local expression \eqref{xsk} and using notation
\eqref{eq:nablaiterated}, we obtain that $X$ is a newtonoid of
order $k$ if and only if for $\alpha \in \{1,...,k\}$ we have
$$h_{\alpha}X=\frac{1}{\alpha !} \nabla^{(\alpha)}X^i
\frac{\delta}{\delta y^{(\alpha)i}},$$ which is equivalent to
formula \eqref{eq:symadaptedbasis}.

For a vector field $X \in J^{k-1}\left(\mathfrak{X}^k_S\right)$
let $Z\in \mathfrak{X}^k_S$ be such that $X=J^{k-1}Z$. Using
formula \eqref{knabla} it follows that \begin{eqnarray*} J(\nabla
X) & =& (J\circ \nabla \circ J^{k-1})Z = (J \circ h_{k-1} \circ
\mathcal{L}_S \circ h_{k-1} \circ J^{k-1})Z \\ & =& \left(J \circ
h_{k-1} \circ \mathcal{L}_S \circ (J^{k-1}-v_k \circ
J^{k-1})\right)Z. \end{eqnarray*} Using $J\circ h_{k-1}=v_k \circ
J$, the fact that for $Z\in \mathfrak{X}^k_S$ we have $(J\circ
\mathcal{L}_S \circ J^{k-1})Z=(1-k)J^{k-1}Z$, and the second
formula in \eqref{jlshvk} we obtain $ J(\nabla X)=v_k(X)$.

Suppose $X=J^{k-1}Z$ for some $Z\in
\mathfrak{X}(T^kM\setminus\{0\})$. Starting with $v_k(X)=J(\nabla
X)$ then the arguments used to prove the first implication show
that $v_kJ^k\mathcal{L}_SZ=0$. Then
\begin{eqnarray*}
J^{k-1}\pi^k_SZ  = J^{k-1}(\operatorname{Id}+ J\mathcal{L}_S)Z=
J^{k-1}Z+ J^k\mathcal{L}_SZ=X,\end{eqnarray*} and $X\in
J^{k-1}\operatorname{Im}\pi^k_S = J^{k-1}\mathfrak{X}^k_S$.
\end{proof}
For $\alpha=1$, formula \eqref{eq:nablaiterated} reduces to
formula \eqref{nablaxik1}, which represents the action of the
dynamical covariant derivative $\nabla$ on the horizontal
components of a vector field. Using formula
\eqref{eq:symadaptedbasis} it follows that $\nabla^{(\alpha)}X^i$,
which are given by formula \eqref{eq:nablaiterated}, represent the
components of a newtonoid vector field of order $k$ with respect
to the basis \eqref{abasis}. Therefore $\nabla^{(\alpha)}X^i$,
behave as the components of a vector field on the base manifold,
for arbitrary components $X^i$ that behave as the components of a
vector field on the base manifold. This observation will allow us
to extend the action of $\nabla^{(\alpha)}$ to arbitrary vector
fields on $T^kM\setminus\{0\}$ as follows. Consider a vector field
$X$ on $T^kM\setminus\{0\}$ expressed in the basis \eqref{abasis}
as follows
\begin{eqnarray} X= \sum_{\beta=0}^k X^{(\beta)i}
\frac{\delta}{\delta y^{(\beta)i}}. \label{aXTKM}
\end{eqnarray} Hence, for each $\beta \in \{0,1,...,k\}$, $X^{(\beta)i}$ behave as the components of a
vector field on the base manifold. Accordingly, for each $\alpha
\in \{1,...,k\}$, the components $\nabla^{(\alpha)} X^{(\beta)i}$
will have the same property. Therefore, we can define the
\emph{$\alpha$-th dynamical covariant derivative} of $X$ by
\begin{eqnarray} \nabla^{(\alpha)}X= \sum_{\beta=0}^k \nabla^{(\alpha)} X^{(\beta)i}
\frac{\delta}{\delta y^{(\beta)i}}, \ \alpha \in \{1,...,k\}.
\label{nalphaxtkm}
\end{eqnarray}
Denote $\nabla^{(0)}=\operatorname{Id}$. For each $\alpha \in
\{1,...,k\}$, we set $\nabla^{(\alpha)}f=S^{\alpha}(f)$, for $f\in
C^{\infty}(T^kM\setminus\{0\})$. Then, $\nabla^{(\alpha)}$
satisfies a Leibniz-type rule of order $\alpha$:
\begin{eqnarray*}\nabla^{(\alpha)}(fX)=\sum_{\beta=0}^{\alpha} {\alpha \choose
\beta} \nabla^{(\beta)}f \cdot
\nabla^{(\alpha-\beta)}X.\end{eqnarray*} Using the Leibniz-type
rule of order $\alpha$, and the requirement that
$\nabla^{(\alpha)}$ commutes with tensor contraction, we can
extend the action of $\nabla^{(\alpha)}$ to arbitrary tensor
fields. For example, consider $\omega$ a $1$-form on
$T^kM\setminus\{0\}$. Then, for an arbitrary vector field $X$ on
$T^kM\setminus\{0\}$, we have
\begin{eqnarray}
\left(\nabla^{(\alpha)}\omega\right)(X)= S^{\alpha}(\omega(X)) -
\sum_{\beta=1}^{\alpha} {\alpha \choose \beta}
\left(\nabla^{(\alpha-\beta)}\omega\right)
\left(\nabla^{(\beta)}X\right). \label{naomega}
\end{eqnarray}
We emphasize that, for the moment, $\nabla^{(\alpha)}$ may not be
the $\alpha$-th iteration of the dynamical covariant derivative
$\nabla$, unless a specific nonlinear connection will be chosen.
In Theorem \ref{kcharnabla} we characterize when
$\nabla^{(\alpha)}=\nabla^\alpha$ for all $\alpha \ge 0$.

\begin{prop} \label{prop:knjn}  Consider the dynamical
covariant derivative $\nabla$  associated to a pair $(S,
(h_{\alpha}, v_k))$. A vector field $X\in
\mathfrak{X}(T^kM\setminus\{0\})$ is a dynamical symmetry if and
only if it is a newtonoid and
\begin{eqnarray}
\frac{1}{k!}\nabla\left(J^k\nabla^{(k)}X\right)+\Phi(X)=0.
\label{knjn}
\end{eqnarray}
\end{prop}
\begin{proof}
A vector field $X\in \mathfrak{X}(T^kM\setminus\{0\})$ is a
dynamical symmetry if and only if $h_{\alpha}[S,X]=0$, for all
$\alpha \in \{0,...,k-1\}$, and $v_k[S,X]=0$. The first set of
conditions $h_{\alpha}[S,X]=0$ for all $\alpha \in \{0,...,k-1\}$
(or the first $k$ conditions) are equivalent with $[S,X]\in
\operatorname{Ker} J$, which is equivalent with $X\in
\mathfrak{X}^k_S$.

For any $X\in \mathfrak{X}(T^kM\setminus\{0\})$ we have
\begin{eqnarray*}
\nabla(v_kX)+\Phi(X)=v_k[S,X].
\end{eqnarray*}
For any for $X\in \mathfrak{X}^k_S$ (with $h_0 X=X^i\delta/\delta
x^i$), formulae \eqref{eq:symadaptedbasis} and \eqref{nalphaxtkm}
imply that
\begin{eqnarray*}
   v_kX =\frac{1}{k!}\nabla^{(k)}X^i\frac{\partial}{\partial y^{(k)i}} = \frac{1}{k!}J^k\nabla^{(k)}X.
\end{eqnarray*}
Using the above observations the result follows.
\end{proof}

From formula \eqref{knabladelta} we see that the dynamical
covariant derivative $\nabla$ has a different action on
distributions $H_{0}, \ldots, H_{k-1}$ than on distribution $V_k$.
The next proposition characterize when these actions coincide. For
example, one sufficient condition is that $\nabla J= 0$. When
$k=1$ we know that the canonical nonlinear connection for a
semispray is uniquely determined by $\nabla J = 0$. Proposition
\ref{kcharnabla1} shows that when $k>1$ there is more freedom, and
$\nabla J = 0$ determines only part of a nonlinear connection. To
uniquely determine a non-linear connection from a semispray one
needs stronger conditions. Such conditions will be given in
Theorem \ref{kcharnabla}.

\begin{prop} \label{kcharnabla1} Consider a semispray of order $k$, $S$,
a nonlinear connection $(h_{\alpha}, v_k)$, and the dynamical
covariant derivative $\nabla$ associated to the pair
$(S,(h_{\alpha}, v_k))$. Then, the following conditions are
equivalent:
\begin{itemize} \item[i)] $\nabla$ restricts to a map $\nabla: J^{k-1}(\mathfrak{X}^k_S) \to J^{k-1}(\mathfrak{X}^k_S)$
that satisfies the Leibniz rule with respect to the $\ast$
product;
\item[ii)] $\nabla J^{\alpha}=0$, for some $\alpha \in
\{1,...,k\}$;
\item[iii)] $(k+1)v_k\circ J^{k-1}= J^{k-1} + \mathcal{L}_SJ\circ
J^{k-1}$; \item[iv)] $N^i_{(1)j}={\partial G^i}/{\partial
y^{(k)j}};$
\item[v)] $\nabla=\nabla^{(1)}$, where $\nabla^{(1)}$ is defined in formula \eqref{nalphaxtkm}.
\end{itemize}
\end{prop}
\begin{proof}
We first prove that $\nabla J^k=0$ if and only if $\nabla
J^{\alpha}=0$, for some $\alpha \in \{1,...,k\}$. In view of first
formula \eqref{nablajk} it follows that if $\nabla J^{\alpha}=0$,
for some $\alpha \in \{1,...,k\}$, then $\nabla J^k=0$. For the
converse we will show that $\nabla J^k=0$ implies $\nabla J=0$,
which will also imply $\nabla J^{\alpha}=0$, for all $\alpha \in
\{1,...,k\}$. From first formula \eqref{nablajk}, we have that
$\nabla J^k=\nabla J \circ J^{k-1}$ and therefore the $(1,1)$-type
tensor field $\nabla J$ vanishes on vector fields in
$\operatorname{Im} J^{k-1}=V_{k-1}$. From local formulae
\eqref{nablajalpha} it follows that $\nabla J$ vanishes on vector
fields in $\bigoplus_{\alpha=0}^{k-2}H_{\alpha}$, and hence
$\nabla J=0$.

Using the above considerations, formulae \eqref{nablajk} and
\eqref{nablajalpha}, and the identity
$\mathcal{L}_SJ^k=\mathcal{L}_SJ\circ J^{k-1} - (k-1)J^{k-1}$ we
have that conditions $ii)$, $iii)$, and $iv)$ are equivalent.

We prove implication $ii) \Longrightarrow i)$. Consider $X\in
J^{k-1}\left(\mathfrak{X}^k_S\right)$. According to Lemma
\ref{lem:kvnablax}, we have that $v_kX=J\nabla X$. If we apply
$\nabla$ to both sides of this equality and use $\nabla v_k=0$ and
$\nabla J=0$ we obtain $v_k(\nabla X)=J\nabla \left(\nabla
X\right)$, which by Lemma \ref{lem:kvnablax} implies that $\nabla
X\in J^{k-1}\left(\mathfrak{X}^k_S\right)$. Hence, $\nabla$
preserves the set  $J^{k-1}(\mathfrak{X}^k_S)$.

For $f\in C^{\infty}(T^kM\setminus\{0\})$ and $X\in
J^{k-1}\left(\mathfrak{X}^k_S\right)$, formula \eqref{kast}
reduces to $f\ast X=fX+S(f)JX$. Using the fact that $\nabla$
satisfies the Leibniz rule for the usual $\cdot$ product it
follows that \begin{eqnarray} \nabla(f\ast X)-\nabla(f)\ast
X-f\ast \nabla X=\nabla(f)\left(\nabla J\right)X. \label{knfast}
\end{eqnarray} Since $\nabla J=0$ it follows that $\nabla$
satisfies the Leibniz rule for the $\ast$ product.

We prove now the implication $i) \Longrightarrow ii)$. By formula
\eqref{nablajalpha}, which gives the local expression for $\nabla
J^{\alpha}$, we only need to show that $\nabla J=0$ on vector
fields in $\operatorname{Im} J^{k-1}$. A set of generators for
vector fields in $\operatorname{Im} J^{k-1}$ is given by
$J^{k-1}\left(\mathfrak{X}^k_S\right)\cup
\left(\mathfrak{X}^{v_k}(T^kM\setminus\{0\})\right)$. Since
$\nabla J$ vanishes on
$\left(\mathfrak{X}^{v_k}(T^kM\setminus\{0\})\right)$, it remains
to show that $\nabla J$ vanishes on
$J^{k-1}\left(\mathfrak{X}^k_S\right)$. Since $\nabla$ satisfies
the Leibniz rule for the $\ast$ product, using formula
\eqref{knfast} it follows that $S(f)\left(\nabla J\right)X=0$, for
an arbitrary function $f\in C^{\infty}(T^kM\setminus\{0\})$ and an
arbitrary vector field $X\in
J^{k-1}\left(\mathfrak{X}^k_S\right)$. Therefore, $\nabla J=0$ on
the set $J^{k-1}\left(\mathfrak{X}^k_S\right)$.

Using formulae \eqref{knabladelta}, \eqref{eq:nablaiterated}, and
\eqref{nalphaxtkm} we have that $\nabla^{(1)}=\nabla$ if and only
if $\nabla J=0$. Hence conditions $ii)$ and $v)$ are equivalent.
\end{proof}
Up to this point, the structure of Section \ref{sec:hode} follows
the structure of Section \ref{sec:sode}, and each result for a
system of SODE has a counterpart for a system of HODE. Due to the
complications imposed by the geometry of a system of HODE we will
have results that are specific only to the higher order case.

Although Proposition \ref{kcharnabla1} fixes only a part of the
nonlinear connection, we can obtain now an explicit covariant form
for the equations of variation \eqref{symmk}, given by the
following Jacobi-type equations \eqref{inv_symmk}.

\begin{prop} \label{kjacobi1}
Consider $\nabla$ the dynamical covariant derivative associated to
a pair $(S,(h_{\alpha}, v_k))$, and suppose that $\nabla J = 0$.
Then, a vector field $X$ on $\mathfrak{X}(T^kM\setminus\{0\})$ is
a dynamical symmetry if and only if $X$ is a newtonoid and
\begin{eqnarray}
\frac{1}{k!}\nabla\left(\nabla^{(k)}X^i\right) +
\sum_{\alpha=0}^{k-1}\frac{1}{\alpha!}R^{i}_{(\alpha)j}\nabla^{(\alpha)}X^j=0,
\label{inv_symmk}
\end{eqnarray}
where $X^i$ are the horizontal components of $X$ such that  $h_0 X
= X^i\delta/\delta x^i$.
\end{prop}

\begin{proof}
If $X$ is a newtonoid, then formulae \eqref{localphi} and
\eqref{eq:symadaptedbasis} imply that
$$
\Phi(X) = \sum_{\alpha=0}^{k-1} \frac{1}{\alpha!}R^{i}_{(\alpha)j}
\nabla^{(\alpha)}(X^j )\frac{\partial}{\partial y^{(k)i}},
$$
and the result follows by Proposition \ref{prop:knjn}.
\end{proof}

Note that while using formula \eqref{knjn} in Proposition
\ref{prop:knjn}, to characterize a newtonoid $X$ on
$T^kM\setminus\{0\}$, we essentially need the assumption $\nabla
J=0$. This implies that the two actions of $\nabla$ in formulae
\eqref{knabladelta} coincide. Therefore, the action of $\nabla$ on
the components of the vertical vector field $v_kX$ is given by the
second formula \eqref{nablaxik1}.

Next theorem will completely determine the nonlinear connection
and this will allow us to provide simpler expressions for the $k$
curvature components, $R^{i}_{(\alpha)j}$, of the Jacobi
endomorphism $\Phi$.

\begin{thm} \label{kcharnabla} Consider a semispray of order $k$, $S$, a nonlinear connection $(h_{\alpha}, v_k)$,
and the dynamical covariant derivative $\nabla$ associated to the
pair $(S,(h_{\alpha}, v_k))$. Then, the following conditions are
equivalent:
\begin{itemize} \item[i)] $\mathcal{L}_SJ+\operatorname{Id}-
(k+1)v_k=i_J\Phi$;
\item[ii)] $\nabla^{(\alpha)}=\nabla^{\alpha}$, for all $\alpha \in \{1,
2,...,k\}$;
\item[iii)] $M_{(1)i}^j=\partial G^j/\partial y^{(k)i}$ and $\alpha M^i_{(\alpha)j}=S\left(M_{(\alpha-1)j}^i\right)
+ M_{(\alpha-1)j}^pM^i_{(1)p}$,  for all $\alpha \in \{2,...,k\}$;
\item[iv)] $N_{(1)i}^j=\partial G^j/\partial y^{(k)i}$ and $\alpha N^i_{(\alpha)j}=S\left(N_{(\alpha-1)j}^i\right)
- N^i_{(\alpha-1)p}N^p_{(1)j}$,  for all $\alpha \in \{2,...,k\}$.
\end{itemize}
\end{thm}
\begin{proof}
Using formulae \eqref{localphi} and \eqref{lsjk} we obtain that
equality i) holds true if and only if $\nabla J=0$ and for each
$\alpha \in \{1,...,k-1\}$ we have
\begin{eqnarray}
R^i_{(\alpha)j}=(k+1)\left(\frac{\delta G^i}{\delta y^{(\alpha)j}}
- N^i_{(k+1-\alpha)j}\right). \label{riaj1} \end{eqnarray} In view
of formulae \eqref{localr}, for $\alpha$ in $\{1,...,k-1\}$, above
equations \eqref{riaj1} are equivalent to
\begin{eqnarray}
(k+1-\alpha)N^i_{(k+1-\alpha)j} = S(N^i_{(k-\alpha)j}) -
N^s_{(1)j} N^i_{(k-\alpha)s} - \sum_{\beta=2}^{k-\alpha}
I^s_{(\beta)j} N^i_{(k+1-\alpha-\beta)}. \label{nsi}
\end{eqnarray}

When $\alpha \in \{1,..., k-1\}$, equality \eqref{nsi} is
equivalent with $I^i_{(k-\alpha+1)j}=0$. Hence condition i) holds
if and only if $\nabla J=0$ and  $I^i_{(\alpha)j}=0$ for all
$\alpha \in \{2,..., k\}$.   Induction and equation \eqref{ialpha}
show that the latter conditions are equivalent with the second
formulae in iv).

Let us first note that condition ii) is equivalent with $\nabla J
=0$ and $\nabla^{(\alpha)} = \nabla \circ \nabla^{(\alpha-1)}$ for
all $\alpha\in \{2, \ldots, k\}$. For $\alpha\in \{1, \ldots,
k\}$, and $\beta\in\{0,...,k\}$, equations
\eqref{eq:nablaiterated} and \eqref{nalphaxtkm} imply that
\begin{eqnarray*}
\nabla^{(\alpha)}\frac{\delta}{\delta
y^{(\beta)j}}=\alpha!M^i_{(\alpha)j} \frac{\delta}{\delta
y^{(\beta)i}}, \label{nalphadelta}
\end{eqnarray*}
and when $\alpha\in \{2, \ldots, k\}$ we have
\begin{eqnarray*}
\left(\nabla \circ \nabla^{(\alpha-1)}\right)\frac{\delta}{\delta
y^{(\beta)j}}=(\alpha-1)!\left(S\left(M_{(\alpha-1)j}^i\right) +
M_{(\alpha-1)j}^pM^i_{(1)p}\right)\frac{\delta}{\delta
y^{(\beta)i}}. \label{nnalphadelta}
\end{eqnarray*}
It follows that ii) and iii) are equivalent.

Using formulae \eqref{nm} and \eqref{naomega}, and induction for
$\alpha \in \{1,...,k\}$, it follows that $\nabla^{(\alpha)}$ has
the following action on the dual basis \eqref{dbasis}
\begin{eqnarray*}
\nabla^{(\alpha)}\delta y^{(\beta)i}=-\alpha! N^i_{(\alpha)j}
\delta y^{(\beta)j}, \quad \beta\in\{0,1,...,k\}.
\label{nalphadual}
\end{eqnarray*}
From this formula, using the action of the dynamical covariant
derivative $\nabla$ on $1$-forms, we obtain, for each
$\beta\in\{0,1,...,k\}$,
\begin{eqnarray*}
\left(\nabla \circ \nabla^{(\alpha-1)}\right)\delta
y^{(\beta)i}=(\alpha-1)!\left(S\left(N_{(\alpha-1)j}^i\right) -
N_{(\alpha-1)p}^iN^p_{(1)j}\right)\delta y^{(\beta)j}.
\end{eqnarray*}
In view of the above two formulae it follows that condition ii)
for $\alpha \in \{2,...,k\}$ and the second formulae in condition
iv) are equivalent.
\end{proof}

The $(1,1)$-type tensor $\mathcal{L}_SJ -
\operatorname{Id}+(k+1)v_k$ in the first item of Theorem
\ref{kcharnabla} can be used to characterize some other nonlinear
connections. For example if we require it to vanish, we obtain a
nonlinear connection studied in \cite{andres91, bucataru05,
cantrijn86, crampin06}.

For a semispray $S$ of order $k$, we say that the nonlinear
connection characterized by any of the four conditions in Theorem
\ref{kcharnabla} is called the \emph{canonical nonlinear
connection} of $S$. For the remaining of this section we will
consider this canonical nonlinear connection induced by a
semispray $S$. The dual coefficients $M^i_{(\alpha)j}$ of this
connection were considered first by Miron and Atanasiu in
\cite{mironat96}. Next theorem characterizes the symmetries of a
semispray $S$ using the dynamical covariant derivative $\nabla$
associated to a semispray $S$ and its canonical nonlinear
connection.

\begin{thm} \label{knablasym} Let $(h_{\alpha}, v_k)$ be the canonical nonlinear connection induced by a semispray $S$,
and let $\nabla$ be the corresponding dynamical covariant
derivative.
\begin{itemize}
\item[i)] A vector field $X$ on $T^kM\setminus\{0\}$ is a dynamical symmetry
if and only if $X$ is a newtonoid and
\begin{eqnarray}
\frac{1}{k!} \nabla^{k+1}J^kX + \Phi(X)=0. \label{ksyminv}
\end{eqnarray} \item[ii)] A vector field $X\in \mathfrak{X}(M)$ is a Lie
symmetry if and only if
\begin{eqnarray}
\frac{1}{k!} \nabla^{k+1}X^{k,0} + \Phi(X^{k,k})=0.
\label{kliesyminv}
\end{eqnarray}
\item[iii)] Locally, a newtonoid vector field $X$ on $T^kM\setminus\{0\}$ is a dynamical
symmetry of $S$ if and only if its horizontal components satisfy
the Jacobi equation
\begin{eqnarray}
\frac{1}{k!}\nabla^{k+1}X^i +
\sum_{\alpha=0}^{k-1}\frac{1}{\alpha!}
R^{i}_{(\alpha)j}\nabla^{\alpha}X^j=0. \label{ksymcoord}
\end{eqnarray}
\item[iv)] The components of the Jacobi endomorphism are given by
\begin{eqnarray}
\label{krij} R^{i}_{(\alpha)j} &=&(k+1)\left(\frac{\delta
G^i}{\delta y^{(\alpha)j}} - N^i_{(k+1-\alpha)j}\right), \quad \alpha \in\{1,...,k-1\} \\
\nonumber R^{i}_{(0)j}&=& (k+1)\frac{\delta G^i}{\delta x^{j}} -
S\left(N^i_{(k)j}\right) + N^i_{(k)l}N^l_{(1)j}.
\end{eqnarray} \end{itemize}
\end{thm}
\begin{proof}
Part i) follows by combining Proposition \ref{prop:knjn} and
Theorem \ref{kcharnabla}. Part ii) follows from Part i) using the
fact that $J^kX^{k,k}=X^{k,0}$. Part iii) represents the local
expression of Part i).  In part iv), the first formula in
\eqref{krij} was proven in the proof of Theorem \ref{kcharnabla}.
The second formula in \eqref{krij} follows by setting $\alpha = 0$
and $I_{(2)j}^i=0,\ldots, I_{(k-1)j}^i=0$ in formula
\eqref{localr}.
\end{proof}

\subsection{First order variation for systems of HODE}
Let $c: I \to M$ be a geodesic of a semispray $S$ of order $k$. A
\emph{geodesic variation} of $c$ is a smooth map $V: I\times
(-\varepsilon, \varepsilon) \to M$, $V=V(t,s)$ such that
\begin{itemize}
\item[i)] $V(t,0)=c(t)$, for all $t$ in $I$,
\item[ii)] $V(t,s)$ is a geodesic for all $s$ in $(-\varepsilon,
\varepsilon)$. \end{itemize}
\begin{defn} \label{def:jacobik}
A vector field $\xi:I \to TM$, along a geodesic $c=\pi \circ \xi$
of a semispray $S$, is called a \emph{Jacobi field} of $S$ if it
is the variation vector field $\xi(t)=\left. \partial_s
V(t,s)\right|_{s=0}$ of a geodesic variation $V$.
\end{defn}
Next we  provide  a sufficient condition for a vector field
$\xi\colon I\to TM$ to be a Jacobi field. This sufficient
condition generalizes the traditional Jacobi equation to a
semispray of order $k$.

For a vector field $\xi:I \to TM$, $\xi(t)=\left(x^i(t),
\xi^i(t)\right)$, along a regular curve $c=\pi \circ \xi$, and
$\alpha \in \{0,...,k\}$, we define the lifted vector fields
$\xi^{k,\alpha}:I \to T(T^kM\setminus\{0\})$, along the $k$-th jet
lift $j^kc$ of $c$, as follows
$$\xi^{k,\alpha}=\left(x^i(t),
\frac{1}{1!}\frac{dx^i}{dt}(t), ...,
\frac{1}{k!}\frac{d^kx^i}{dt^k}(t), 0,...,0, \xi^i(t),
\frac{1}{1!} \frac{d\xi^i}{dt}(t), ...,\frac{1}{\alpha!}
\frac{d^{\alpha}\xi^i}{dt^{\alpha}}(t) \right).$$ Since $\xi$ is a
vector field along $c$, for each $t\in I$, we can find a vector
field $X\in \mathfrak{X}(M)$ such that $X\circ c=\xi$ near $t$.
Then $\xi^{k,\alpha} = X^{k,\alpha}\circ j^kc$, for all $\alpha
\in \{0,...,k\}$. We assume now that $c$ is a geodesic, which
means that $j^kc$ is an integral curve of a semispray $S$. For
$\alpha\ge 0$ we define
\begin{eqnarray}
\nabla^\alpha \xi^{k,0}=\left(\nabla^\alpha X^{k,0}\right)\circ
j^kc. \label{nablaxik}
\end{eqnarray}
If locally $\xi=(x^i, \xi^i)$ and $X=X^i {\partial}/{\partial
x^i}$, then $\nabla \xi^{k,0}= \nabla \xi^i
\partial/\partial y^{(k)i}$, where
\begin{eqnarray}
\nabla \xi^i = \left(\nabla X^i\right)\circ j^kc =
\frac{d\xi^i}{dt}+ N_{(1)j}^i\left(x,
\frac{1}{1!}\frac{dx}{dt},..., \frac{1}{k!}
\frac{d^kx}{dt^k}\right)\xi^j. \label{nablaxik0}
\end{eqnarray}
By induction on $\alpha$ we find that
\begin{eqnarray*}
\nabla^{\alpha} \xi^i = \left(\nabla^{\alpha} X^i\right)\circ
j^kc, \quad \alpha\ge 0,
\end{eqnarray*}
where $\nabla^{\alpha} X^i$ are iterations of formula
\eqref{nablaxik1}, that correspond to formulae
\eqref{eq:nablaiterated} for the nonlinear connection determined
by Theorem \ref{kcharnabla}. Hence, using $\nabla^\alpha X^{k,0} =
\nabla^\alpha X^i {\partial}/{\partial y^{(k)i}}$,
\begin{eqnarray*}
   \nabla^{\alpha}\xi^{k,0}&=&\nabla^\alpha \xi^i \frac{\partial}{\partial y^{(k)i}}, \quad \alpha\ge 0.
\end{eqnarray*}
It follows that definition \eqref{nablaxik} does not depend on the
choice of $X$, and components $\nabla^{\alpha} \xi^i$ transform as
a tensor on $M$.

\begin{prop} \label{prop:invjacobik}
Let $\xi$ be a vector field $\xi:I \to TM$ along a geodesic of a
semispray $S$ of order $k$. If $\xi$ is a \emph{Jacobi field} of
$S$ then it satisfies the following \emph{Jacobi equation}
\begin{eqnarray}
\frac{1}{k!} \nabla^{k+1}\xi^{k,0} + \Phi(\xi^{k,k})=0.
\label{jacobixik0}
\end{eqnarray}
If locally $\xi(t)=\left(x^i(t), \xi^i(t)\right)$, equation
\eqref{jacobixik0} is equivalent to
\begin{eqnarray}
\frac{1}{k!}\nabla^{k+1}\xi^i +
\sum_{\alpha=0}^{k-1}\frac{1}{\alpha!}
R^{i}_{(\alpha)j}\left(x,\frac{1}{1!}\frac{dx}{dt}, ...,
\frac{1}{k!}\frac{d^kx}{dt^k} \right)\nabla^{\alpha}\xi^j=0.
\label{eq:invjacobik}
\end{eqnarray}
\end{prop}

\begin{proof}
By induction on $\alpha$ we have that
\begin{eqnarray*}
\delta y^{(\alpha)i} \xi^{k,k} &=&
\frac{1}{\alpha!}\nabla^{\alpha} \xi^i, \quad\alpha\in
\{0,\ldots,k\}.
\end{eqnarray*}
Hence formulae \eqref{jacobixik0} and \eqref{eq:invjacobik} are
equivalent. To prove that formula \eqref{jacobixik0} holds, let
$c\colon I \to M$ be the geodesic $c=\pi\circ \xi$ of $S$, let
$t\in I$, and let $X$ be a vector field  $X\in \mathfrak{X}(M)$
such that $X\circ c=\xi$ on some interval $I_0\ni t$. Since $\xi$
is a Jacobi field there is a geodesic variation $V(t,s)$ such that
$V(t,0)=c(t)$, ${\partial V^i}/{\partial s}(t,0)=\xi^i(t),$ and
$$
\frac{\partial^{k+1} V^i}{\partial t^{k+1}}+ (k+1)!G^i\left(V,
\frac{1}{1!}\frac{\partial V}{\partial t},...,
\frac{1}{k!}\frac{\partial^k V}{\partial t^k}\right)=0, \forall
(t,s)\in I_0 \times (-\varepsilon, \varepsilon).
$$
Differentiating the last equation with respect to $s$ and setting
$s=0$ gives
\begin{eqnarray*}
\frac{d^{k+1}\xi^{i}}{dt^{k+1}} + (k+1)! \left(\frac{\partial
G^{i}}{\partial x^{j}} \xi^{j} + \frac{\partial G^{i}}{\partial
y^{(1)j}} \frac{1}{1!}\frac{d\xi^{j}}{dt} + \cdots +
\frac{\partial G^{i}}{\partial y^{(k)j}} \frac{1}{k!}
\frac{d^{k}\xi^{j}}{dt^{k}}\right) &=& 0.
\end{eqnarray*}
Since $c$ is a geodesic, it follows that $S^\alpha (f)\circ j^k(c)
= {d^\alpha (f\circ j^kc)}/{dt^\alpha}$ for all $\alpha \ge 0$ and
$f \in C^\infty(T^kM\setminus\{0\})$. Thus
\begin{eqnarray*}
   \left(S^{k+1}(X^i) + (k+1)!X^{k,k}(G^i)\right)\circ j^kc =  0.
\end{eqnarray*}
By Lemma \ref{lem33} and inclusion \eqref{inclusionk} it follows
that $J[S,X^{k,k}]=0$. Then
$$
  v_k[S, X^{k,k}] = \frac{1}{k!} S^{k+1}(X^i) + X^{k,k} (G^i),
$$
and $[S,X^{k,k}]\circ j^kc=0$. Repeating the argument in
Proposition \ref{prop:knjn} shows that
$$
  \frac{1}{k!} \nabla^{k+1} (J^k X^{k,k} )\circ j^kc + \Phi(X^{k,k}) \circ j^kc=0.
$$
Equation \eqref{jacobixik0} follows.
\end{proof}
A proof for the converse of Proposition \ref{prop:invjacobik} will
require to extend all the techniques developed in \cite{bd10}.
Another aspect that will have to be addressed in the future is the
role of the curvature components $R^i_{(\alpha)j}$ for the
geodesic behavior.

\section{Applications}
\label{sec:applic} In this section we motivate the applicability
of the geometric theory developed in this paper using examples
from various fields such as: the equivalence problem, the inverse
problem of the calculus of variations, and biharmonicity.

First, we consider a second order Lagrangian, $L_2$, \cite{bm09}
derived from a Riemannian structure. The corresponding
Euler-Lagrange equations form a system of fourth order
differential equations, whose solutions are biharmonic curves
\cite{cadeo06}. For this system of differential equations, the
components of the Jacobi endomorphism are functions of the
curvature of the Riemannian metric and its dynamical covariant
derivatives. This will motivate the use of \emph{curvature
components} for the components of the Jacobi endomorphism.

We also show that the curvature components of the Jacobi
endomorphism are useful to express geometric invariants that were
associated previously to third or fourth order ordinary
differential equations. See \cite{crampin06, dridi06, fels96,
neut02}.

\subsection{Prolongation of a Riemannian structure} \label{subsec:prol}
In this section we start with a Riemannian space $(M, g)$ and
construct a $2$-nd order Lagrangian $L_2\in C^{\infty}(T^2M)$
\cite{bm09}. The Euler-Lagrange equations for $L_2$ determine a
system of fourth order ordinary differential equations and hence a
semispray of order $3$. Solutions of this system are biharmonic
curves for the Riemannian space \cite{cadeo06}. The $3$ components
$R^i_{(\alpha)j}$, $\alpha \in \{0,1,2\}$, of the Jacobi
endomorphism can be expressed in terms of the curvature components
$R^{i}_{jkl}$ and their first and second order dynamical covariant
derivatives.

Consider $g=g_{ij}(x)dx^i\otimes dx^j$ a semi-Riemannian metric on
$M$ and denote by $\gamma^i_{jk}(x)$ its Christoffel symbols. The
geodesic spray of the Riemannian metric $g$ has  the following
components of the Jacobi endomorphism $R^i_j(x,
y^{(1)})=R^{i}_{kjl}(x)y^{(1)k}y^{(1)l}$. We denote
\begin{eqnarray}
z^{(2)i}= y^{(2)i}+ \frac{1}{2}\gamma^i_{jk}(x)y^{(1)j}y^{(1)k}.
\label{z2}
\end{eqnarray}
It follows that $z^{(2)i}$ behave as the components of a vector
field on $M$. These components were interpreted as the covariant
form of acceleration in \cite[(6.5)]{bm09} and as half of the
components of the tension field in \cite{cadeo06}. Therefore, the
function $L_2: T^2M \to \mathbb{R}$, given by
\begin{eqnarray}
L_2(x,y^{(1)}, y^{(2)})=\frac{1}{2}g_{ij}z^{(2)i}z^{(2)j},
\label{l2}
\end{eqnarray} is a second order Lagrangian, \cite{bm09, cadeo06}. The variational
problem for $L_2$ leads to the following Euler-Lagrange equations
\begin{eqnarray}
\frac{\partial L_2}{\partial x^i} -
\frac{d}{dt}\left(\frac{\partial L_2}{\partial y^{(1)i}}\right) +
\frac{1}{2}\frac{d^2}{dt^2} \left(\frac{\partial L_2}{\partial
y^{(2)i}}\right)=0, \label{el2}
\end{eqnarray}
which form a system of fourth order ordinary differential
equations. Consider $S\in \mathfrak{X}(T^3M)$
\begin{eqnarray}
S=y^{(1)i}\frac{\partial}{\partial x^i} +
2y^{(2)i}\frac{\partial}{\partial y^{(1)i}} +
3y^{(3)i}\frac{\partial}{\partial y^{(2)i}} -
4G^i\frac{\partial}{\partial y^{(3)i}}, \label{s3} \end{eqnarray}
the corresponding semispray of order $3$. This semispray is
uniquely determined by the following equation
\begin{eqnarray}
\frac{\partial L_2}{\partial x^i} - S\left(\frac{\partial
L_2}{\partial y^{(1)i}}\right) + \frac{1}{2}S^2
\left(\frac{\partial L_2}{\partial y^{(2)i}}\right)=0.
\label{els3}
\end{eqnarray}
Indeed, the functions $G^i$, of the semispray $S$, uniquely
determined by equation \eqref{els3}, are given by
\begin{eqnarray}
6g_{ji}G^j=\frac{\partial L_2}{\partial x^i} -
d_T\left(\frac{\partial L_2}{\partial y^{(1)i}}\right) +
\frac{1}{2}d_T^2 \left(\frac{\partial L_2}{\partial
y^{(2)i}}\right), \label{gi3} \end{eqnarray} where $d_T$ is the
Tulczyjew operator on $T^3M$. Although, above expression for the
coefficients $G^i$ of the semispray $S$ is not so easy to handle
with, we will show that one can obtain the corresponding nonlinear
connection, which is determined by any of the four equivalent
conditions in Theorem \ref{kcharnabla}. In formula \eqref{els3} we
apply $\partial/\partial y^{(3)k}$, use the fact that the
Lagrangian $L_2$ and its partial derivatives are functions on
$T^2M$ we obtain
$$\frac{\partial G^i}{\partial y^{(3)j}} = N^i_{(1)j}(x,y^{(1)}) =
M^i_{(1)j}(x,y^{(1)}) = \gamma^i_{kj}(x)y^{(1)k}.$$ For the other
two coefficients $N^i_{(2)j}(x, y^{(1)}, y^{(2)})$, $N^i_{(3)j}(x,
y^{(1)}, y^{(2)}, y^{(3)})$, of the nonlinear connection we use
formulae iv) from Theorem \ref{kcharnabla}.

Consider $\nabla$ the dynamical covariant derivative associated to
the semispray $S$ and the corresponding nonlinear connection. It
follows that
$$z^{(2)i}=\frac{1}{2}\nabla y^{(1)i},$$ hence the name of
\emph{covariant form of acceleration} for $z^{(2)i}$. Equation
\eqref{els3} is then equivalent to
\begin{eqnarray}
\frac{\delta L_2}{\delta x^i} - \nabla\left(\frac{\delta
L_2}{\delta y^{(1)i}}\right) + \frac{1}{2}\nabla^2
\left(\frac{\delta L_2}{\delta y^{(2)i}}\right)=0.
\label{nablaels3}
\end{eqnarray}
Using the following formulae
\begin{eqnarray*}
\frac{\delta z^{(2)i}}{\delta y^{(1)j}}= 0, \ 2g^{ik}\frac{\delta
L}{\delta x^k}=R^i_{k}z^{(2)k}, \ 2\left( \frac{\delta
z^{(2)i}}{\delta x^j} + \gamma^i_{jk}z^{(2)k}\right)= R^i_j,
\end{eqnarray*} we can rewrite equation \eqref{nablaels3} as
follows
\begin{eqnarray} \nabla^2 z^{(2)i} + R^i_j z^{(2)j}=0.\label{n2z2}
\end{eqnarray} Note that this equation is equivalent to
$$  \nabla^3 y^{(1)i} + R^i_j \nabla y^{(1)j}=0.$$
Therefore the Euler-Lagrange equations \eqref{el2} can be written
as follows
$$  \nabla^3 \left(\frac{dx^i}{dt}\right) + R^i_j \nabla \left(\frac{dx^j}{dt}\right)=0,$$
which is the system of differential equations for biharmonic
curves \cite[(1)]{cadeo06}. Since a geodesic of the Riemannian
metric $g$ is a solution of the system $\nabla(dx^i/dt)=0$, it
follows that any geodesic is a biharmonic curve. Conversely, if
the sectional curvature is non-positive then any biharmonic closed
curve is a geodesic, see \cite{jiang86}. However, if we change the
parameter of an arbitrary geodesic by a third order polynomial,
the new curve is no longer a geodesic but a biharmonic curve, see
\cite{balmus09}.

Using any of the equivalent form \eqref{n2z2}, \eqref{gi3}, or
\eqref{nablaels3}, of equations \eqref{els3}, we can express the
components $R^i_{(\alpha)j}$, $\alpha \in\{0,1,2\}$ of the Jacobi
endomorphism in terms of the curvature tensor $R^i_{jkl}$ of the
Riemannian metric. Two of them are
\begin{eqnarray*} 3R^i_{(2)j} & = & 2 R^i_j= 2R^i_{sjk}y^{(1)s}y^{(1)k}, \\
3R^i_{(1)j} & = & \nabla R^i_j + R^i_{jsk} \nabla y^{(1)s}
y^{(1)k},
\end{eqnarray*} The other component $R^i_{(0)j}$, has a more complicated formula, but it is a function of the curvature tensor
$R^i_{jkl}$, its first, and second order dynamical derivative.

\subsection{The Wuenschmann invariant}
\label{subsec:Wuenschmann}

A geometric invariant that relates third order ordinary
differential equations with certain classes of conformal Lorentz
metrics on three dimensional manifolds, was proposed by K.
Wuenschmann in his Ph.D thesis of 1905. The role of the
Wuenschmann invariant for the geometry of a third order
differential equation under contact transformations was discussed
in \cite{neut02}. See also \cite{dridi06} for a generalization of
the problem to fourth order ordinary differential equations, where
the equivalence problem leads to three invariants. A geometric
method for obtaining the Wuenschmann invariant for a system of
third order ordinary differential equations was proposed in
\cite{crampin06}. In this section we show that the Wuenschmann
invariant can be obtained from the components of the Jacobi
endomorphism.

On a $1$-dimensional manifold $M$, consider a third-order
differential equation
$$\frac{d^3x}{dt^3}+3!G\left(x, \frac{1}{1!}\frac{dx}{dt},
\frac{1}{2!}\frac{d^2x}{dt^2}\right)=0,$$ and the corresponding
semispray $S$ of order $2$. There exists an associated conformal
Lorentzian structure on the $3$-dimensional solution space of this
equation if and only if an invariant vanishes, \cite{crampin06,
frittelli06}. This invariant is called the \emph{Wuenschmann
invariant} and it is given by \cite[eq. (4)]{neut02}
\begin{eqnarray*}
W_3 = -\frac{1}{2}S^2\left(\frac{\partial G}{\partial
y^{(2)}}\right) - 3\frac{\partial G}{\partial y^{(2)}}
S\left(\frac{\partial G}{\partial y^{(2)}}\right) + 3
S\left(\frac{\partial G}{\partial y^{(1)}}\right) - 2
\left(\frac{\partial G}{\partial y^{(2)}}\right)^3  + 6
\frac{\partial G}{\partial y^{(1)}}\frac{\partial G}{\partial
y^{(2)}} - 6\frac{\partial G}{\partial x}. \end{eqnarray*}

Consider the two curvature components given by formula
\eqref{krij} of the Jacobi endomorphism:
\begin{eqnarray*}
R_{(0)} & = & 3\frac{\partial G}{\partial x}- 3\frac{\partial
G}{\partial y^{(1)}}\frac{\partial G}{\partial y^{(2)}}
-\frac{1}{2}S^2\left(\frac{\partial G}{\partial y^{(2)}}\right) +
\left(\frac{\partial G}{\partial y^{(2)}}\right)^3, \\
R_{(1)} & = & 3\frac{\partial G}{\partial y^{(1)}}-
\frac{3}{2}\left(\frac{\partial G}{\partial y^{(2)}}\right)^2  -
\frac{3}{2}S\left(\frac{\partial G}{\partial y^{(2)}}\right).
\label{r0r1}
\end{eqnarray*}
It follows that the Wuenschmann invariant and the components of
the Jacobi endomorphism are related as follows \begin{eqnarray}
W_3=\nabla R_{(1)}- 2R_{(0)}. \label{w3r} \end{eqnarray} It is
important to note that same formula \eqref{w3r} was obtained in
\cite{crampin06} for two other different nonlinear connections.

\subsection{The inverse problem of Lagrangian mechanics for scalar fourth-order
ordinary differential equations} \label{subsec:inv}

In \cite{fels96}, Fels has shown that a scalar fourth-order
ordinary differential equation admits a variational multiplier if
and only if two invariants vanish. These two geometric invariants
were associated with a fourth-order equation using Cartan's
equivalence method. One of these invariants appears also in the
list of invariants proposed by Dridi and Neut in \cite{dridi06}
for studying the equivalence problem for fourth order differential
equations under fiber preserving diffeomorphisms. In this section
we show that this invariant can be expressed in terms of the
components of the Jacobi endomorphism using a very similar formula
as in the previous section.

On a $1$-dimensional manifold $M$, consider a fourth-order
differential equation
$$\frac{d^4x}{dt^4}+4!G\left(x, \frac{1}{1!}\frac{dx}{dt},
\frac{1}{2!}\frac{d^2x}{dt^2}, \frac{1}{3!}\frac{d^3x}{dt^3}
\right)=0,$$ and the corresponding semispray $S$ of order $3$.
Consider the following Wuenschmann-type invariant
\begin{eqnarray*}
W_4 = -\frac{2}{3}S^2\left(\frac{\partial G}{\partial
y^{(3)}}\right) - 4\frac{\partial G}{\partial y^{(3)}}
S\left(\frac{\partial G}{\partial y^{(3)}}\right) + 4
S\left(\frac{\partial G}{\partial y^{(2)}}\right) - \frac{8}{3}
\left(\frac{\partial G}{\partial y^{(3)}}\right)^3 + 8
\frac{\partial G}{\partial y^{(2)}}\frac{\partial G}{\partial
y^{(3)}} - 8\frac{\partial G}{\partial y^{(1)}}. \end{eqnarray*}
The invariant $W_4$ is related to the invariant $I_1$ considered
by Fels in \cite{fels96} as follows $W_4=-3I_1$. Also $(8/3)W_4$
is the numerator of $I_9$ in \cite{dridi06}.

Consider two curvature components $R_{(1)}$ and $R_{(2)}$ given by
formula \eqref{krij} of the Jacobi endomorphism:
\begin{eqnarray*}
R_{(1)} & = & 4\frac{\partial G}{\partial y^{(1)}}-
4\frac{\partial G}{\partial y^{(2)}}\frac{\partial G}{\partial
y^{(3)}} - \frac{2}{3}S^2\left(\frac{\partial G}{\partial
y^{(3)}}\right) + \frac{4}{3}
\left(\frac{\partial G}{\partial y^{(3)}}\right)^3, \\
R_{(2)} & = & 4\frac{\partial G}{\partial y^{(2)}}- 2
\left(\frac{\partial G}{\partial y^{(3)}}\right)^2 - 2
S\left(\frac{\partial G}{\partial y^{(3)}}\right). \label{r1r2}
\end{eqnarray*}
It follows that the invariant $W_4$ and the components of the
Jacobi endomorphism are related as follows \begin{eqnarray}
W_4=\nabla R_{(2)}- 2R_{(1)}. \label{w4r} \end{eqnarray} We note
the similarities for expressing the two invariants $W_3$ and $W_4$
in formulae \eqref{w3r} and \eqref{w4r}.

As an example consider the classical spinning particle
\cite{deLeon92}. The motion of a particle rotating about a
translating center is governed by the following system of fourth
order differential equations
\begin{eqnarray}
\frac{d^{4}x^{i}}{dt^{4}}+\omega^{2}\frac{d^{2}x^{i}}{dt^{2}}=0,\quad
i\in\{1,2,3\}, \label{spinp}\end{eqnarray} where $\omega$ is a
real, non-zero constant. For this system of fourth order we have
\begin{eqnarray*} G^{i}=\frac{1}{12}\omega^{2}y^{(2)i}, \quad
N_{(1)j}^{i}=N_{(2)j}^{i}=N_{(3)j}^{i}=0. \end{eqnarray*} The only
non-vanishing component of the Jacobi endomorphism is
\begin{eqnarray*}
R_{(2)j}^{i}=\frac{1}{3}\omega^{2}\delta_{j}^{i}.\end{eqnarray*}
Therefore, the Wuenschmann invariant is $W_{4}=\nabla R_{(2)}-
2R_{(1)}=0$. Note that the system \eqref{spinp} is separable and
hence we can view each of the three equations of the system
separately. For each of them, the Wuenschmann invariant vanishes
and therefore each equation is a Lagrangian equations. A second
order Lagrangian for the system \eqref{spinp} is the following
\begin{eqnarray*} L(x,y^{(1)},
y^{(2)})=\delta_{ij} y^{(2)i}y^{(2)j} -
\frac{\omega^2}{2}\delta_{ij} y^{(1)i}y^{(1)j}.\end{eqnarray*} For
the system \eqref{spinp}, the Jacobi equations are
\begin{eqnarray*}
\nabla^{4}\xi^{i}+\omega^{2}\nabla^{2}\xi^{i}=0,\end{eqnarray*}
which can be integrated. Therefore we can obtain all geodesic
variations of the system \eqref{spinp}.

\begin{acknowledgement*}
This work has been supported by: the Romanian Ministry of
Education (grant CNCSIS - UEFISCSU, PNII - IDEI 398), the Academy
of Finland (project 13132527 and Centre of Excellence in Inverse
Problems Research), and by the Institute of Mathematics at Aalto
University.
\end{acknowledgement*}


\begin{thebibliography}{99}

\bibitem{andres91} Andres, L.C., de Le\'on, M., and Rodriguez, M.: \emph{Connections
on tangent bundles of higher order associated to regular
Lagrangians}. Geometriae Dedicata, \textbf{39} (1991), 17--28.

\bibitem{antonelli03} Antonelli, P.L., Bucataru, I.: \emph{ KCC-theory of a System of Second
Order Differential Equations}. Handbook of Finsler Geometry, vol.
I, Kluwer Academic Publisher, (2003), 83--174.

\bibitem{balmus09} Balmu\c s, A: \emph{Biharmonic maps and submanifolds}, Differential
Geometry - Dynamical Systems Monographs, \textbf{10}, Geometry
Balkan Press, Bucharest, 2009.

\bibitem{bao00} Bao, D., Chern, S.-S., Shen, Z.: \emph{An introduction to
Riemann-Finsler geometry}, Springer, 2000.

\bibitem{bucataru00} Bucataru, I.: \emph{Horizontal lifts in the higher
order geometry}. Publicationes Mathematicae, \textbf{56} (1-2)
(2000), 21--32.

\bibitem{bucataru05} Bucataru, I.: \emph{Linear connections for systems of
higher order differential equations}. Houston Journal of
Mathematics, \textbf{31} (2) (2005), 315--332.

\bibitem{bucataru07a} Bucataru, I.: \emph{Metric nonlinear connections}.
Differential Geometry and its Application, \textbf{35} (3) (2007),
335--343.

\bibitem{bucataru07b} Bucataru, I.: \emph{Canonical semispray for higher order
Lagrange spaces}. Comptes Rendus Math\'ematique. Acad\'emie des
Sciences. Paris, \textbf{345} (2007), 269--272.

\bibitem{bm09} Bucataru, I., Miron, R.: \emph{The geometry of systems of third order
differential equations induced by second order Lagrangians}.
Mediterranean Journal of Mathematics, \textbf{6} (4) (2009),
483--500.

\bibitem{bd09} Bucataru, I., Dahl, M.F.: \emph{ Semi-basic 1-forms and Helmholtz
conditions for the inverse problem of the calculus of variations}.
Journal of Geometric Mechanics, \textbf{1} (2) (2009), 159--180.

\bibitem{bd10} Bucataru, I., Dahl, M.F.: \emph{A complete lift for
semisprays}. International Journal of Geometric Methods in Modern
Physics, \textbf{7} (2) (2010), 267--287.

\bibitem{byrnes96} Byrnes, G. B.: \emph{A linear connection for higher-order ordinary
differential equations}.  Journal of Physics A, \textbf{29} (8)
(1996), 1685--1694.

\bibitem{cadeo06} Caddeo, R., Montaldo, S., Oniciuc, C., Piu, P.:
\emph{The Euler-Lagrange method for biharmonic curves}. Mediterr.
J. Math., \textbf{3} (3-4) (2006),  449--465.

\bibitem{carinena91} Cari\~nena, J.F., Mart\'inez, E.: \emph{Generalized Jacobi equation and inverse problem in
classical mechanics}, in ``Group Theoretical Methods in Physics"
(eds. V. V. Dodonov and V. I. Manko), Proc. 18th Int. Colloquim
1990, Moskow, vol. II, Nova Science Publishers, (1991) New York.

\bibitem{docarmo92} do Carmo, M. P.: \emph{Riemannian geometry},
Birkh\"auser Boston, 1992.

\bibitem{catz74} Catz, G., \emph{Gerbes et connexions sur le fibre tangent
d'ordre 2}, C. R. Acad. Sc. Paris, \textbf{278} (1974), 347--349.

\bibitem{crampin71} Crampin, M.: \emph{On horizontal distributions on the
tangent bundle of a differentiable manifold}. J. London Math.
Soc., \textbf{2} (3) (1971), 178--182.

\bibitem{cantrijn86} Crampin, M., Sarlet, W., Cantrijn, F.: \emph{Higher
order differential equations and higher order Lagrangian
Mechanics}. Proc. Camb. Phil. Soc., \textbf{99} (1986), 565--587.

\bibitem{crampin96} Crampin, M., Mart\'inez, E. Sarlet, W.: \emph{Linear connections for
systems of second-order ordinary differential equations}.  Ann.
Inst. H. Poincar\'e Phys. Th\'eor., \textbf{65} (2) (1996),
223--249.

\bibitem{crampin06} Crampin, M., Saunders, D. J.: \emph{On the geometry of higher-order
ordinary differential equations and the Wuenschmann invariant}.
Groups, geometry and physics, Monogr. Real Acad. Ci. Exact.
Fís.-Quím. Nat. Zaragoza, 29, Acad. Cienc. Exact. Fís. Quím. Nat.
Zaragoza, Zaragoza, (2006), 79--92.

\bibitem{dridi06} Dridi, R., Neut, S.: \emph{The equivalence problem for fourth order differential
equations under fiber preserving diffeomorphisms}. Journal of
Mathematical Physics, \textbf{47} (2006), 013501.

\bibitem{fels96} Fels, M. E.: \emph{The inverse problem of the calculus of variations
for scalar fourth-order ordinary differential equations}. Trans.
Amer. Math. Soc. 348 (1996), 5007--5029.

\bibitem{frittelli06} Frittelli, S., Kozameh, C., Newman, E. T.,
Nurowski, P.; \emph{Differential equations and Cartan
connections}. Topics in mathematical physics, general relativity
and cosmology in honor of Jerzy Pleban\~nski, World Sci. Publ.,
Hackensack, NJ, (2006), 193--200.

\bibitem{grifone72} Grifone, J.: \emph{Structure presque-tangente et connexions I}.
Ann. Inst. Fourier, \textbf{22} (1) (1972), 287--334.

\bibitem{grifone00} Grifone, J., Muzsnay, Z.: \emph{Variational Principles for Second-order Differential Equations},
World-Scientific, 2000.

\bibitem{jiang86} Jiang, G.Y.: \emph{2-harmonic maps and their first and second
variational formulas}, Chinese Ann. Math. Ser. A, \textbf{7}
(1986), 389--402.

\bibitem{kosambi35} Kosambi D.D.: \emph{Systems of Differential Equations of
the Second Order}. Quart. J. Math. Oxford Ser., \textbf{6} (21)
(1935), 1--12.

\bibitem{kosambi36} Kosambi D.D.: \emph{Path Spaces of higher order}. Quart. J. Math.
Oxford Ser., \textbf{7} (26) (1936), 97--104.

\bibitem{deleon95} de Le\'on, M., de Diego, D. M.: \emph{Symmetries and constants of the motion for higher-order
Lagrangian systems}. Journal of Mathematical Physics, \textbf{36}
(8)(1995), 4138--4161.

\bibitem{deleon85} de Le\'on, M., Rodrigues, P. R.: \emph{Generalized classical
mechanics and field theory. A geometrical approach of Lagrangian
and Hamiltonian formalisms involving higher order derivatives}.
North-Holland Mathematics Studies, 112. Notes on Pure Mathematics,
102. North-Holland Publishing Co., Amsterdam, 1985.

\bibitem{deLeon92} de Le\'on, M., Rodrigues, P. R.: \emph{The inverse problem
of Lagrangian dynamics for higher-order differential equations: a
geometrical approach}. Inverse Problems, \textbf{8} (4) (1992),
525--540.

\bibitem{marmo86} Marmo, G., Mukunda, N.: \emph{Symmetries and constants of the
motion in the lagrangian formalism on TQ: beyond point
transformations}. Nuovo Cim. B, \textbf{92} (1986) 1--12.

\bibitem{miron94} Miron, R.: \emph{Noether theorem in higher order
Lagrangian mechanics}. International Journal of Theoretical
Physics, \textbf{34} (7) (1994), 1123--1146.

\bibitem{miron97} Miron, R.: \emph{The geometry of higher-order
Lagrange spaces. Applications to mechanics and physics}. Kluwer
Academic Publishers, 1997.

\bibitem{mirona94} Miron, R., Anastasiei, M.: \emph{The Geometry of Lagrange Spaces: Theory and
Applications}. Kluwer Academic Publishers,  1994.

\bibitem{mironat96} Miron, R., Atanasiu, Gh.: \emph{Prolongation of Riemannian,
Finslerian and Lagrangian structures}. Rev. Roumaine Math. Pures
Appl., \textbf{41} (3-4) (1996), 237--249.

\bibitem{nicolaecsu07} Nicolaescu, L.: \emph{Geometry of
Manifolds}, World Scitific, 2007.

\bibitem{neut02} Neut, S., Petitot, M.: \emph{La g\'eom\'etrie de l'\'equation
$y'''=f(x,y',y'')$}. C. R. Acad. Sci. Paris, Ser. I, \textbf{335}
(2002), 515--518.

\bibitem{punzi09} Punzi, R., Wohlfarth, M.R.: \emph{Geometry and stability of
dynamical systems}. Physical Review E, \textbf{79} (2009), 046606.

\bibitem{sarlet90} Sarlet, W.: \emph{Adjoint symmetries of second-order differential equations and
generalizations}.  Differential geometry and its applications
(Brno, 1989) World Sci. Publ., (1990)  412--421.

\bibitem{saunders02} Saunders, D.J.: \emph{On the inverse problem
for even-order ordinary differential equations in the higher-order
calculus of variations}. Differential geometry and its
Applications, \textbf{16} (2002), 149--166.

\bibitem{shen01} Shen, Z.: \emph{Differential geometry of spray and
Finsler spaces}, Springer, 2001.

\bibitem{szilasi03} Szilasi, J.: \emph{A setting for spray and Finsler
geometry}, Handbook of Finsler geometry. Kluwer Acad. Publ.,
Dordrecht, Vol. 2, 2003, 1183--1426.

\bibitem{tulczyjew76} Tulczyjew, W.M.: \emph{The Lagrange differential
geometry}. Bull. Acad. Polon. Sci. \textbf{24} (1976), 1089--1096.

\bibitem{yano73} Yano, K., Ishihara, S.: \emph{Tangent and cotangent bundles}, Marcel Dekker, Inc., 1973.


\end{thebibliography}
\end{document}